\def\versiondate{22 Apr. 2020}
\input math.macros

\let\nobibtex = t
\let\noarrow = t
\input eplain
\beginpackages
\usepackage{url}
   \usepackage{color}  
\endpackages

\input Ref.macros

\checkdefinedreferencetrue
\continuousfigurenumberingtrue
\theoremcountingtrue
\sectionnumberstrue
\forwardreferencetrue
\citationgenerationtrue
\nobracketcittrue
\hyperstrue
\initialeqmacro

\hyphenation{Kai-man-ov-ich}

\input kv.key
\bibsty{../../texstuff/myapalike}

\def\bfz{{\bf 0}}

\def\F{\scr F}

\def\gp{\Gamma}

\def\invar{{\cal I}}
\def\sinvar {{\cal J}} 
\def\dist{{\rm dist}}
\def\lamps{\phi}
\def\Lamps{\Phi}


\def\gpe{\gamma}
\def\ent{{\bfi h}}

\def\lincr{\Psi}
\def\lit{\mathop {\rm lit}}
\def\rad{\mathop {\rm rad}}
\def\lampg{\frak{L}}
\def\baseg{\Lambda}
\def\id{{\ss id}}
\def\bfid{{\bf ID}}
\def\oid{{\bf o}}
\def\enum#1{{\bf #1}}
\def\ranen#1{{\cal #1}}
\def\FF{{\Bbb F}}

\def\MM{{\cal M}} 

\def\BH{{\bf BH}}
\def\bnd{{\bf b}}
\def\chngo#1{\delta^{#1}}  

\def\grnmet{\rho}  
\def\iprod#1#2{\langle #1, #2\rangle}
\def\Cov{\mathop {\rm Cov}}
\def\ewset{Q}
\def\verts{\gp}
\def\bp{o}

\def\TT{{\Bbb T}}
\def\cut{{\ss cut}}
\def\tmax{t_{\rm max}}
\def\pd{\pi}  
\def\walk{\Xi}
\def\mub{\mu_{\rm base}}  
\def\sgp{\Delta}  
\def\as{{\ss Z}}   
\def\dmax{\dist_{\max}}  
\def\rl#1{#1}

\ifproofmode \relax \else\head{}
{Version of \versiondate}\fi
\vglue20pt

\title{Poisson Boundaries of Lamplighter Groups:}
\smallskip
\title{Proof of the Kaimanovich--Vershik Conjecture}

\author{Russell Lyons and Yuval Peres}

\abstract{We answer positively a question of Kaimanovich and Vershik
from 1979, showing that the final configuration of lamps for simple random
walk on the lamplighter group over $\Z^d$ ($d \ge 3$) is the Poisson
boundary. For $d \ge 5$, this had been shown earlier by Erschler (2011). We
extend this to walks of more general types on more general groups.}

\bottomIII{Primary
20F69, 
60B15, 
60J50. 
Secondary
43A05, 
20F65. 
}
{Random walks,
free metabelian group, entropy, harmonic functions.}
{Research partially supported by NSF grants DMS-1007244 and DMS-1612363
and Microsoft Research.}

\bsection{Introduction}{s.intro}

Suppose that $\gp$ is a countable infinite
group and $\mu$ is a probability measure
on $\gp$ whose support generates $\gp$ (as a group).
A function $f \colon \gp \to \R$ is called \dfn{harmonic} if $f(x) =
\sum_z \mu(z) f(x z)$ for all $x \in \gp$.
If all bounded harmonic functions are
constant, then $(\gp, \mu)$ is said to have the \dfn{Liouville property}.
A general theory for the non-Liouville case was initiated by
\refbmulti{Furst:Annals,Furst:ICM,Furst:Lie}, who defined the notion of
Poisson boundary to describe the set of bounded harmonic functions.
Such harmonic functions are closely linked to the \dfn{$\mu$-walk}, which is
the Markov chain with transition probabilities $p(x, y) := \mu(x^{-1} y)$.
Earlier work on boundaries for general Markov chains is due to \ref
b.Blackwell/, \ref b.Feller/, \ref b.Doob/, \ref b.Hunt/, and \ref
b.Feldman/; a special case for groups was established by \ref b.DynMal/.

\ref b.Rosenblatt/ and \ref b.KV/ \rl{(announced in \ref b.KV:announce/)}
proved a conjecture of \ref b.Furst:boundary/ that
$\gp$ is amenable iff there is a symmetric $\mu$ whose support generates $\gp$
such that $(\gp, \mu)$ is Liouville.
Another open question had been whether there exists an amenable group
with a symmetric non-Liouville measure.
To answer this, \refbmulti{KV:announce,KV}
utilized certain restricted wreath products $\Z_2 \wr \Z^d$, now commonly
called lamplighter groups, where $\Z_2 := \Z/(2\Z)$ is referred to as the
lamp group and $\Z^d$ as the base group.
These are solvable (hence amenable)
groups of exponential growth.
To define them more generally, let $\lampg$ and $\baseg$ be two groups.
Then $\lampg \wr \baseg$ is the semidirect product $\Big(\sum_{z \in \baseg}
\lampg\Big) \rtimes \baseg$, where $\baseg$ acts on $\sum_{z \in \baseg} \lampg$ by
$$
(x \Lamps)(z) := \Lamps(x^{-1} z)
\,.
$$
Thus, if $(\Lamps, x), (\lincr, y) \in \sum_{x \in \baseg} \lampg \times
\baseg$,
then
$$
(\Lamps, x) (\lincr, y)
=
\big(\Lamps \cdot (x \lincr), xy\big)
\,.
$$
The interpretation of an element $(\Lamps, x)$ is
that a lamplighter is at $x$, there is one lamp at each element of
$\baseg$, each lamp has a state in $\lampg$, and
$\Lamps$ gives the
states of all the lamps.
If $\lampg$ and $\baseg$ are both finitely generated, then so is their
restricted wreath product.
To see this,
write $o$ for the identity in $\baseg$ and $\id$ for the identity in $\lampg$.
Write $\bfid$ for the function
that is equal to $\id$ identically on $\baseg$.
Also, write $\chngo s$ for the element of $\sum_{z \in \baseg} \lampg$ that
equals $s$ at $o$ and equals $\id$ elsewhere; thus, $\bfid = \chngo \id$.
If $S_1$ and $S_2$ are generating sets for $\lampg$ and
$\baseg$, respectively, then an often-used generating set for $\lampg \wr
\baseg$
is $\big\{(\chngo {s_1}, o) \st s_1 \in S_1\big\} \cup \big\{(\bfid, s_2) \st
s_2 \in S_2\big\}$.
Multiplying $(\Lamps, x)$ on the right by a generator $(\chngo {s_1}, o)$
changes the state of the lamp at $x$ by $s_1$, while multiplying $(\Lamps,
x)$ on the right by a generator $(\bfid, s_2)$ moves the lamplighter to $x
s_2$.
Since every element of $\sum_{z \in \baseg} \lampg$
is the identity of $\lampg$ at all but finitely many $z \in \baseg$, the
above set does indeed generate $\lampg \wr \baseg$.

Let $\mu$ be a finitely supported, symmetric probability measure whose
support generates $\Z_2 \wr \Z^d$.
\ref b.KV/, Proposition 6.4, showed that
$(\Z_2 \wr \Z^d, \mu)$ is Liouville iff $d \le 2$.
\refbmulti{KV:announce,KV}
also asked for a description of the
Poisson boundary for finitely supported $\mu$ on the lamplighter groups
$\Z_2 \wr \Z^d$ when it is nontrivial, \rl{which, in the symmetric case,
amounts to} $d \ge 3$.
Moreover, they suggested a natural candidate, namely, $(\Z_2)^{\Z^d}$ with
the probability measure given by the final configuration of lamps under the
$\mu$-walk.
On $\Z_2 \wr \Z^d$,
the final configuration of lamps,
which we will denote by $\Lamps_\infty$,
exists because the projection of the walk to the
base group $\Z^d$ is transient.

In 2008, a breakthrough was achieved by
\refbmulti{Erschler,Erschler:ICM}, who
proved that the conjecture of \ref b.KV:announce/ is correct
when $d \ge 5$.

We show here that the conjecture of \ref b.KV:announce/ is correct for all
$d \ge 3$.
In fact, we prove the following \rl{main result. Say that a probability
measure $\mu$ on $\lampg \wr \baseg$ has \dfn{bounded lamp range} if $\{x
\in \baseg \st \texists {(\Lamps, y) \in \lampg \wr \baseg} \mu(\Lamps, y)
> 0 \hbox{ and } \Lamps(x) \ne \id\}$ is a finite set. This means that one
step of the random walk can change the lamp values on only a set of bounded
size, which holds, for example, if $\mu$ has finite support.
Write $\mub$ for the projection of $\mu$ on $\baseg$.
\procl t.finiteent
Let $\lampg$ be a nontrivial finite group and
$\baseg$ be a finitely generated, infinite group.
Let $\mu$ be a probability measure on $\lampg \wr \baseg$
with finite entropy and bounded lamp range,
and whose support generates $\lampg \wr \baseg$.
If $\mub$ generates a transient random walk on $\baseg$, 
then the Poisson boundary is $\lampg^\baseg$ endowed with the law of\/
$\Lamps_\infty$.
\endprocl
It follows readily from known results that if the projected measure $\mub$
generates a recurrent random walk on $\baseg$ and $(\lampg, \nu)$ is
Liouville for every $\nu$ whose support generates $\lampg$,
then the Poisson boundary of
$(\lampg \wr \baseg, \mu)$ is trivial whenever the support of $\mu$
generates $\lampg \wr \baseg$; see \ref p.recurrent/.
This was proved earlier for abelian lamp groups
(Proposition 1.2 in \ref B.Kaim:examples/) and, more generally, nilpotent
lamp groups (Theorem 3.1 of \ref B.Kaim:solvable/).
\vskip\parskip
\ref t.finiteent/
is proved in \ref s.jumps/, with various minor strengthenings. 
Theorems \briefref t.not-semigp/ and \briefref t.finiteent-nonLiou/
give settings in which the assumptions that $\lampg$ be finite
and that $\baseg$ be finitely generated can be removed.
The assumptions that $\lampg$ is finite and
$\mu$ has bounded lamp range are replaced by a
second-moment assumption in \ref t.2ndmoment/ when $\baseg = \Z^d$.
}

Entropy is a key quantity in the study of Poisson boundaries.
\rl{We are aware of no significant results that identify a
nontrivial Poisson boundary in the presence of infinite entropy, although
\ref b.ForTio/ manage to reduce finite logarithmic moment to finite
entropy on free semigroups.}

%



We introduce an enhanced version of the celebrated entropy criterion
of \ref b.Kaim:AOM/, which has been the key tool for identification of
Poisson boundaries. This is presented in \ref c.vadim-quant2/ and used in
\ref s.SRW/. We also discuss it informally below in the context of the
history of the subject.

Poisson boundaries are related to other important aspects of random walks.
One fundamental aspect is to determine, given a random walk on a group
$\gp$,
its set of possible asymptotic behaviors, by which we mean the
$\sigma$-field $\invar$ on the path space $\gp^\N$ invariant under time shifts.
There is a well-known correspondence between $\invar$ and the space $\BH$
of bounded harmonic functions on $\gp$.
In particular, the invariant $\sigma$-field is trivial (i.e., consists only
of sets of probability 0 or 1) iff all bounded harmonic functions are
constant.

\rl{Following} the introduction of asymptotic entropy by
\refbmulti{Avez:ent,Avez:Choquet,Avez:croiss,Avez:harmonic}
and \rl{the 0--2 law of}
\ref b.Derriennic:0-2/,
a foundational
paper by \ref b.KV/, announced in \ref b.KV:announce/, developed a general
theory to analyze Poisson boundaries.
In particular, Avez, \ref b.Derriennic:applic/, and Kaimanovich--Vershik
proved that if
$\mu$ has finite entropy, then the Avez (asymptotic) entropy of the
$\mu$-walk is 0 iff the walk is Liouville.
\ref b.Varopoulos/ showed that for finitely supported, symmetric
$\mu$, the rate of escape of the $\mu$-walk is sublinear iff $(\gp, \mu)$
is Liouville. This was extended by \ref b.KarlLed/ to symmetric $\mu$ with
finite first moment with respect to the word metric for a finite generating
set.

\ref b.Erschler:Liouville/ showed \rl{(1)}
that every finitely generated solvable group
of exponential growth admits a symmetric non-Liouville measure, and 
\rl{(2)} that every non-degenerate measure on $\lampg \wr \baseg$ whose
projection to $\baseg$ is transient has nonzero Avez entropy. She also
proved a result similar to \rl{(2)}
for the free metabelian
groups $\FF_d/\FF''_d$ with $d \ge 3$.
Furthermore,
\ref b.Erschler:subexp/ showed that there are groups of
intermediate growth with finite-entropy, symmetric, non-Liouville measures.
\ref b.FHTVF/ extended this to show
that every finitely generated group that is not of
polynomial growth admits a finite-entropy symmetric non-Liouville measure.

\ref b.Furst:Lie/ and \ref b.KV/ gave entropy criteria for identifying the
Poisson boundary.
Two notable papers by \refbmulti{Ledrappier:relation,Ledrappier:discrete}
used this criterion to determine the Poisson boundary for discrete matrix
groups.
\ref b.LedBall/ developed further the entropy method in the context of
rank-one manifolds.
\refbmulti{Kaim:max,Kaim:hyper,Kaim:AOM} refined the
entropy method more generally
and introduced a powerful general
criterion for equality of a given boundary
and the Poisson boundary \rl{(see \ref t.Kcrit-full/ and
\ref c.Kcrit-epsilon/ here)}.

Informally, Kaimanovich's
criterion says that in order that a candidate boundary be
the Poisson boundary, it suffices to find a sequence of random finite sets
$Q_n \subset \gp$, that depend on points of the candidate boundary, such that
$|Q_n|^{1/n} \to 1$ as $n \to\infty$ and $\P[\walk_n \in Q_n]$ is bounded
below, where $\walk_n$ is the location of the random walk at time
$n$.
One of Kaimanovich's important observations was that the sets
$Q_n$ can often be defined geometrically.
This led to his well-known strip and ray criteria.
We enhance Kaimanovich's more
general criterion so that it suffices that $\walk_m \in Q_n$ for {\it
some} $m \ge n$.

Kaimanovich's criteria
led to much progress in identifying Poisson boundaries, such as the
works by
\refbmulti{Kaim:mapping,KaimMasur:teich},
\ref b.Karlsson:Floyd/,
\ref b.Malyutin:locfree/,
\ref b.KarlssonWoess:trees/,
\ref b.Sava/, 
\ref b.BroffSchap/,  
\ref b.GautMath/, 
\ref b.MR3644015/, 
\ref b.NevoSageev/, 
\ref b.MaherTiozzo/, 
and
\ref b.MalSvet/. 

Using these methods, \ref b.Kaim:survey/ made some progress on the
lamplighter question by showing that
for $\mu$ whose projection on the base group, $\Z^d$, has nonzero mean, the
final lamps do indeed give the Poisson boundary.
This problem of identifying the Poisson boundary
has been raised repeatedly (e.g., \ref b.Kaim:solvable/, \ref
b.Vershik:survey/, \ref b.KarlssonWoess:trees/, \ref b.Sava/,
\refbmulti{Erschler,Erschler:ICM}, \ref b.Agelos/) and has been
considered a major open problem in the field.

Beyond \refbauthor{Erschler}'s result on the Kaimanovich--Vershik conjecture,
similar results have been
established for random walks $\mu$ of
finite first moment whose support generates one of
the following groups $\lampg \wr \baseg$:
\beginbulletitems
\bullitem \rl{$\lampg$ is finitely generated and nontrivial,}
$\baseg$ has subexponential growth,
and there is a homomorphism $\psi \colon \baseg \to \Z$ such that if $\pi
\colon \lampg \wr \baseg \to \baseg$ is the canonical projection, then
$(\psi\pi)_* \mu$ has nonzero mean (\ref b.Kaim:survey/);
\bullitem  $\lampg$ is finite and $\baseg$ is
a group with a Cayley graph being a tree of degree at least 3 (\ref
b.KarlssonWoess:trees/);
\bullitem  $\lampg = \Z_2$ and $\baseg$ \rl{is finitely generated and} has
infinitely many ends or is non-elementary hyperbolic (\ref b.Sava/).
\endbulletitems
\noindent
In all these cases, the projection of the random walk to $\baseg$ has linear
rate of escape,
and this makes the analysis considerably simpler.

\ref b.Erschler/ also extended her result on $\lampg \wr \Z^d$ ($\lampg$
finitely generated
and $d \ge 5$) beyond finitely supported $\mu$ to those with finite
third moment, and noted that similar techniques work for free metabelian
groups $\FF_d/\FF_d''$ when $d \ge 5$.

Prior to the work of
\ref b.Erschler/, Kaimanovich's entropy criterion was used in a mostly
geometric fashion that did not require detailed knowledge of the
probabilistic behavior of the random walks.
\refbauthor{Erschler} succeeded in her results by discovering how to leverage
such knowledge of random walks in $\Z^d$ for $d \ge 5$.
In particular, she relied heavily on the existence of a positive
density of cutpoints
(for simple random walk---and analogous behavior in general).
That is, for the lamplighter random walk $\Seq{\widehat X_n}$ on $\Z_2 \wr
\Z^d$, its projection $X_n$ at time $n$ to the base $\Z^d$ is a cutpoint with
probability bounded below over all $n$.  This allowed \refbauthor{Erschler}
to define the required random finite sets $Q_n$ that capture $X_n$ with
probability bounded below.

Our enhanced criterion allows the use of cut-spheres, which do not occur
with positive density, but they do occur infinitely often for $d = 3, 4$.
Use of cut-spheres also simplifies considerably the definition of the
random sets $Q_n$. This is a general feature of our enhanced criterion,
which we illustrate with a simple proof of a conjecture of \ref b.Sava:thesis/.
However, we do not use our enhanced criterion to handle general base
groups, where other innovations are used.
The innovation that is most closely related to cutpoints is to use upper
bounds on the Green function in order to bound the
number of times at which the future of the walk can get close to the
locations of the past of the walk.
In order to handle more general base groups beyond $\Z^d$, other
innovations convert small entropy growth found in various places
into enumeration with small exponential growth of the required sets $Q_n$.

We begin with the definition of the Poisson boundary and
\refbauthor{Kaim:AOM}'s criterion in \ref s.lemmas/.
In order to present the proof of the original conjecture of \ref
b.KV:announce/ in the briefest manner,
we prove \ref t.finiteent/
in \ref s.SRW/
in the special case where the $\mu$-walk
is simple random walk on $\lampg \wr \Z^d$.
We then prove the full \ref t.finiteent/ in \ref s.jumps/.
As did \refbauthor{Erschler}, we consider other step distributions $\mu$ on
$\lampg \wr \Z^d$;
in \ref s.gen/,
we extend her result to $d \ge 3$
and to $\mu$ having finite second moment.
In this broader setting where generators can change lamps arbitrarily far
from the location of the lamplighter,
some technical condition is needed to ensure existence of the limiting lamp
configuration, as discussed at the end of \ref s.gen/.
\rl{One can also ask about infinitely generated base groups, $\baseg$; some
of our results apply in that case: see Theorems \briefref t.not-semigp/ and
\briefref t.finiteent-nonLiou/.}
In \ref s.meta/, we give some details about metabelian groups and similar
groups and discuss our extensions to them.


\bsection{Preliminaries}{s.lemmas}

For a discrete probability distribution $\pd$ on a set $S$, write $H(\pd) :=
\rl{-}\sum_{s \in S} \pd(s) \log \pd(s)$ for the \dfn{entropy} of $\pd$.
For a $\sigma$-field $\F$ and a discrete random variable $X$, write $H(X)$
for the entropy of the distribution of
$X$ and $H(X \mid \F)$ for the conditional
entropy of $X$ given $\F$:
$$
H(X \mid \F):=
-\Ebigg{\sum _{x} \P[X=x \mid \F] \cdot \log \P[X=x \mid \F] } \,.
$$
Our Markov chains will begin at a fixed point; when that point is $x$, we
use $\Psub_x$ for the corresponding probability measure. Usually $x$ will be
the identity element, $\oid$, of a group, $\gp$.
When a transition matrix is given,
we often regard $\Psub_x$ as the law on $\gp^\N$ of the trajectory
of the corresponding Markov chain $\Seq{\walk_n \st n \ge 0}$.
The $\sigma$-field of shift-invariant events is denoted by
$\invar$.
We say that two $\sigma$-fields are equal \dfn{mod 0} if their completions are
equal, generally with respect to $\Psub_\oid$.
The diagonal action of $\gp$ by multiplication
on $\gp^\N$ induces an action of $\gp$ on
$\invar$; a subset $\sinvar \subseteq \invar$ is said to be
\dfn{$\gp$-closed} if $\gpe(A) \in \sinvar$ for all $\gpe \in \gp$ and all
$A \in \sinvar$.

The following criteria of \ref b.Kaim:AOM/ (see Theorem 4.6 and
Corollary 4.6 there, or see Theorem 14.35 and Corollary 14.36 of \ref
b.LP:book/) are essential in identifying Poisson boundaries:

\procl t.Kcrit-full
Let $\Seq{\walk_n \st n \ge 0}$ be a random walk on $\gp$ with
$H(\walk_1) < \infty$.
Let $\invar$ be the associated
invariant $\sigma$-field and $\sinvar \subseteq \invar$ be a
$\gp$-closed sub-$\sigma$-field.
Then $\ent^{\sinvar} :=
\lim_{n \to\infty} n^{-1} H^{\sinvar}(\walk_n)$ converges a.s.\ and in
$L^1$ to the constant $H(\walk_1 \mid \sinvar) - H(\walk_1 \mid \invar)$.
Furthermore, $\ent^{\sinvar} = 0$ iff $\sinvar = \invar$ mod 0.
\endprocl

\procl c.Kcrit-epsilon
Let $\Seq{\walk_n \st n \ge 0}$ be a random walk on $\gp$ with
$H(\walk_1) < \infty$.
Let $\invar$ be the associated
invariant $\sigma$-field and $\sinvar \subseteq \invar$ be a
$\gp$-closed sub-$\sigma$-field.
Suppose that for each $\epsilon > 0$, there is a random sequence
$\Seq{Q_{n, \epsilon} \st n \ge 0}$ of finite subsets of\/ $\gp$ such that
\beginitems
\itemrm{(i)} $Q_{n, \epsilon}$ is $\sinvar$-measurable;
\itemrm{(ii)} $\limsup_{n \to\infty} {1 \over n} \log |Q_{n, \epsilon}| <
\epsilon$ a.s.;
\itemrm{(iii)} $\limsup_{n \to\infty} \Psub_\oid[\walk_n \in Q_{n, \epsilon}] > 0$.
\enditems
Then $\sinvar = \invar$ mod 0.
\endprocl

When $\gp$ is replaced by the lamplighter group $\lampg \wr \baseg$,
we will apply this to the $\lampg \wr \baseg$-closed
$\sigma$-field $\sinvar := \sigma(\Lamps_\infty) \subseteq \invar$ defined
by the limiting configuration of lamps.
Thus, $Q_{n, \epsilon}$ will be a measurable function of configurations
$\lamps_\infty \in \lampg^\baseg$.

In \ref s.SRW/,
we will illustrate the use of a more flexible version of the preceding
corollary, to wit:

\procl c.vadim-quant2
Let $\Seq{\walk_n \st n \ge 0}$ be a random walk on $\gp$ with
$H(\walk_1) < \infty$.
Let $\invar$ be the associated
invariant $\sigma$-field and $\sinvar \subseteq \invar$ be a
$\gp$-closed sub-$\sigma$-field.
Suppose that for each $\epsilon > 0$, there is a random sequence
$\Seq{\rl{\ewset}_{n, \epsilon} \st n \ge 0}$ of finite subsets of\/ $\gp$ such that
\beginitems
\itemrm{(i)} $\ewset_{n, \epsilon}$ is $\sinvar$-measurable;
\itemrm{(ii)} $\limsup_{n \to\infty} {1 \over n} \log |\ewset_{n, \epsilon}| <
\epsilon$ a.s.;
\itemrm{(iii)} $\limsup_{n \to\infty} \Psub_\oid[\texists{m \ge n}   \walk_m
\in \ewset_{n, \epsilon}] > 0$.
\enditems
Then $\sinvar = \invar$ mod 0.
\endprocl

\proof
Write $p_n^{\sinvar}(x, y) := \Psub_x[\walk_n = y \mid \sinvar]$ for the
transition probabilities of the Markov chain conditioned on $\sinvar$.
We will use the following result of \ref b.Kaim:AOM/:
$$
\lim_n {1 \over n} \log p_n^{\sinvar}(\oid,\walk_n) =-\ent^{\sinvar}
\quad \hbox{a.s.}
\label e.condshannon
$$

It suffices to show that $\ent^{\sinvar}=0$. Suppose that
$\ent^{\sinvar}>0$ and define the random $\sinvar$-measurable sets
$$
S_m:=\big\{x \in \verts \st  p_m^{\sinvar}(\oid,x) \le \exp(-m\ent^{\sinvar}/2)
\big\} \,.
$$
For $\epsilon>0$,
$$
\Psubbig_\oid{\walk_m \in \ewset_{n, \epsilon} \cap S_m \bigm| \sinvar} \le
|\ewset_{n, \epsilon}|\cdot \exp(-m\ent^{\sinvar}/2)\, .
$$
Summing over $m \ge n$, we deduce that  for $0<\epsilon<\ent^{\sinvar}/2$,
$$
\Psubbig_\oid{\texists{m \ge n}   \walk_m \in \ewset_{n, \epsilon} \cap S_m
\bigm|
\sinvar} \le   |\ewset_{n, \epsilon}|\cdot c \exp(-n\ent^{\sinvar}/2)
\to 0 \ \hbox{a.s.}
\label e.wset
$$
as $n \to \infty$,  where $c=c(\ent^{\sinvar})$ is a constant. Therefore,
$$
\Psubbig_\oid{\texists{m \ge n}  \walk_m \in \ewset_{n, \epsilon} \cap S_m}
\to 0 \hbox{\quad as } n \to \infty \, .
\label e.wset2
$$
By \ref e.condshannon/, $\Psubbig_\oid{\texists{m \ge n}   \walk_m \notin
S_m} \to 0$ as $n \to \infty$. In conjunction with \ref e.wset2/,
this implies that
$$
\Psubbig_\oid{\texists{m \ge n}   \walk_m \in \ewset_{n, \epsilon}   }   \to 0
\hbox{\quad as } n \to \infty \, ,
$$
contradicting the hypothesis (iii).
\Qed

A \dfn{Poisson boundary} for a random walk on $\gp$ is a quadruple
$(\Theta, \F, \nu, \bnd)$, where $(\Theta, \F, \nu)$ is a probability space
with $\F$ being countably generated and separating points, and where $\bnd
\colon (\gp^\N, \invar) \to (\Theta, \F)$ is a
$\gp$-equivariant measurable map that pushes forward $\Psub_\oid$ to $\nu$ and
such that $\bnd^{-1} \F = \invar$
mod $\Psub_\oid$.
It is unique up to isomorphism.
For more details and background, see
\ref b.Kaim:AOM/ or Definition 14.28 and Theorem 14.29 of \ref b.LP:book/.

When we consider random walks on $\lampg \wr \baseg$, we will write
$\widehat X_n := \walk_n$ and $\widehat X_n =: (\Lamps_n, X_n)$.
Similarly, write $\Seq{\widehat Y_n \st n \ge 1}$ for the increments of
$\Seq{\walk_n \st n \ge 0}$ on $\lampg \wr \baseg$, i.e., $\widehat Y_n
:= \widehat X_{n-1}^{-1} \widehat X_{n}$.
Write $\widehat Y_n =: (\lincr_n, Y_n)$. 
Thus, $\Seq{Y_n \st n \ge 1}$ are IID elements of $\baseg$, used as
increments of the random walk $\Seq{X_n \st n \ge
0}$, i.e., $Y_n := X_{n-1}^{-1} X_{n}$.
Note that while $\Seq{\lincr_n}$ are IID, $\lincr_n$ and $Y_n$ are in general
dependent for each $n$.
Also, for $x \in \baseg$,
$$
\Lamps_n(x)
=
\Lamps_{n-1}(x) \lincr_n(X_{n-1}^{-1} x)
\,.
$$
We generally assume that the support of (the law of) $Y_1$ generates
$\baseg$ and, likewise, the support of $\widehat Y_1$ generates $\lampg \wr
\baseg$.

Let $\lit \lamps$ denote the set of ``lit lamps", $\{x \in \baseg \st
\lamps(x) \ne \id\}$, of $\lamps \in \lampg^\baseg$\rl{, also sometimes
referred to as the support of $\lamps$}.

Suppose that $\Lamps_\infty := \lim_{n \to\infty} \Lamps_n$ exists a.s.
For example, this occurs if $\Ebig{|\!\lit \lincr_1|} < \infty$
and $\Seq{X_n}$ is transient (\ref b.Kaim:solvable/, Theorem 3.3, or
\ref b.Erschler/, proof of Lemma 1.1).
In various cases, we will show that 
\rl{$(\lampg^\baseg, \F, \nu, \bnd)$}
is a Poisson boundary, where
$\F$ is the product $\sigma$-field, $\nu$ is the $\Psub_\bfid$-law of
$\Lamps_\infty$, and
\rl{$\bnd \colon \big((\lampg \wr \baseg)^\N,
\invar, \Psub_\bfid\big) \to (\lampg^\baseg, \F, \nu)$}
takes a sequence to its limiting configuration of lamps; on the set
of measure 0 where the limiting configuration does not exist, we define
$\bnd$ to take the value $\bfid$ for convenience.

We will use $c$ to stand for a positive constant, whose
value can vary from one use to another.

When a group is finitely generated, we use the word metric to define
$|x|$ as the distance between $x$ and the identity element.

\bsection{Proof for the Classical Case}{s.SRW}

Here we give a very short proof of the basic conjecture of \ref
b.KV:announce/ concerning random walks on $\lampg \wr \Z^d$ for $d \ge 3$
and $\lampg$ any nontrivial \rl{finite or countable} group.

\procl t.classical
Let $\lampg$ be a nontrivial \rl{finite or countable} group.
Let $d \ge 3$.
Let $\mu$ be a probability measure of finite entropy on $\lampg \wr \Z^d$
whose support generates $\lampg \wr \Z^d$.
Suppose that $\mu$
is concentrated on
$\big\{(\chngo {s}, o) \st s \in \lampg\big\} \cup \big\{(\bfid, x) \st
x \in \Z^d\big\}$.
If the projection of $\mu$ on $\Z^d$ is finitely supported and has mean $\bfz$,
then the Poisson boundary of $(\lampg \wr \Z^d, \mu)$ is $\lampg^{\Z^d}$
endowed with the law of\/ $\Lamps_\infty$.
\endprocl

\proof
\rl{Write $R$ for the maximum
distance in $\Z^d$ from the current location that one step of the Markov
chain can move.}
For $r>1$, consider the events
$$
\cut_r:=
\big[\texists {m \ge 1} (\all{k<m} |X_k|<r\
\hbox {\rm and} \  \all{j>m} |X_j|>r) \big]
\,.
$$
In the proof of their Proposition 2.1,
\ref b.JamesPeres/ showed that when the projection of $\mu$ is symmetric,
$\Psub_o(\cut_r) \ge c/r$ and $\Psub_o(\cut_r \cap \cut_{r+j}) \le
c/(rj)$. In fact, their proof depends only on estimates of the Green
function, and those hold as long as the projection of $\mu$ has mean
$\bfz$: see, e.g., \ref b.LawLim:book/, Theorem 4.3.1.
Thus, the preceding inequalities of \ref b.JamesPeres/ hold not only for
symmetric $\mu$, but also for those $\mu$ whose projection has mean $\bfz$.
The second moment method applied to $\sum_{r=\rl{R}n}^{n^2} \I{\cut_r}$,
just as in the proof of Proposition 2.1 of \ref b.JamesPeres/,
then yields that
$\Psub_\bp\bigl(\bigcup_{r=Rn}^{n^2} \cut_r \bigr) \ge c (\log n)^2/(\log n)^2 =
c > 0$ \rl{for $n > R$}.
Define
$
\ewset_n :=
\ewset_{n, \epsilon}(\Lamps_\infty)$ to be the set of
$(\lamps,x)$ such that $|x| \le n^2$ and
$$
\lamps(z)= \cases{\Lamps_\infty(z) & if $|z|<|x|$,\cr
                  \id & if $|z|\ge |x|$.\cr
                  }
$$
If $r \le n^2$ and the time $m$ witnesses the event $\cut_r$, then
$\widehat X_m \in \ewset_n$ \rl{and $m \ge r/R$}. 
Therefore,
$\Psubbig_\bfid{\texists{m \ge n} \widehat X_m \in \ewset_n } \ge
\Psub_\bp\bigl(\bigcup_{r=Rn}^{n^2} \cut_r \bigr) \ge c >0$ \rl{for $n >
R$}; since
$|\ewset_n| \le c n^{2d}$, \ref c.vadim-quant2/ implies that
$\sigma(\Lamps_\infty)$ coincides with $\invar$ mod 0.
\Qed

\rl{It is not too hard to extend the above proof to all $\mu$ with finite
support. We leave this as an exercise to the reader who wishes to better
understand the method. A full proof of a more general result is given for
\ref t.2ndmoment/.}

As a further illustration of the usefulness of \ref c.vadim-quant2/, we
prove a conjecture of \ref b.Sava:thesis/.
First we remark that the notion of Poisson boundary extends to all Markov
chains, and criteria such as
\ref c.vadim-quant2/ extend to the setting of
transitive Markov chains: see \ref b.KaimWoess:homog/ for the required
analogues of \ref t.Kcrit-full/ and Equation \ref e.condshannon/, or see
\ref b.LP:book/, Proposition 14.34 and Theorem 14.35.

Now consider the $d$-regular tree, $\TT_d$, and fix an end $\xi$ of
$\TT_d$. The group of graph automorphisms that preserve $\xi$ is known as
the \dfn{affine group of $\TT_d$}; it acts transitively on the vertex set,
$V(\TT_d)$.
Fix some vertex $\bp \in V(\TT_d)$. There is a horodistance
function $d_\xi \colon V(\TT_d) \to \Z$ defined by $d_\xi(\bp) = 0$ and
$d_\xi(x) = d_\xi(y) + 1$ when $y$ is the parent of $x$ (the unique
neighbor of $x$ in the direction of $\xi$).
The affine group preserves differences of values of the horodistance
function.

Let $\lampg$ be a \rl{nontrivial} finite group.
We consider Markov chains $\Seq{\widehat X_n \st n \ge 1} =
\Seq{(\Lamps_n,X_n) \st n \ge 1}$ on the state space
$$
\lampg \wr {\TT_d} :=
\big\{(\lamps, x) \st \lamps \in \lampg^{V(\TT_d)},\
|\!\lit \lamps| < \infty,\  x \in V(\TT_d) \big\}
$$
that change lamps only in a bounded neighborhood of the current location,
make only bounded jumps in the base $\TT_d$, and
whose transition probabilities are invariant under the diagonal action of
the affine group.
Write $R$ for the maximum distance in $\TT_d$ from the current location
that one step of the Markov
chain can move or at which one step of the Markov chain can change the lamps.

\ref b.Sava:thesis/ conjectured the following \ref t.sava/.
She proved that it holds
when $\Ebig{d_\xi(X_1)} \ne 0$ (indeed, with $R < \infty$ replaced by a
first moment condition) or
when $\Seq{X_n}$ is a nearest-neighbor random walk that can change lamps
only at the location of the lamplighter.

\procl t.sava
Let $\Seq{\widehat X_n}$ be a Markov chain that is invariant under the
affine group of\/ $\TT_d$ such that $R < \infty$ and the random
walk projected to the base $\TT_d$ is not constant.
Then the Poisson boundary of $\Seq{\widehat X_n}$ is
$\lampg^{V(\TT_d)}$ endowed with the law of\/ $\Lamps_\infty$.
\endprocl

\proof
We may assume that $\Ebig{d_\xi(X_1)} = 0$. \ref b.CKW/ proved
that $\Seq{X_n}$ converges to $\xi$ a.s.
Let $\xi_n$ be the $\xi$-ancestor of $\bp$ with $d_\xi(\xi_n) = -n$.
Define the cone $C_n := \{x \st \xi_n \hbox{ is an ancestor of } x\}$.

The case when $\Seq{X_n}$ is a nearest-neighbor random walk
is somewhat simpler for our method:
To see how it follows from \ref c.vadim-quant2/, let $\ewset_{n,
\epsilon}(\lamps_\infty)$ be the singleton $\big\{(\lamps_n, \xi_n)\big\}$,
where $\lamps_n(y) = \lamps_\infty(y)$ for $y \in C_n$ and
$\lamps_n(y) = \id$ otherwise.
Let $\alpha := \Psub_x[\all {j \ge 1} X_j \ne x]$; this does not depend on
$x$ by transitivity and is positive by transience.
With $\P_\bfid$-probability 1,
there will be some random smallest time $m \ge n$ such that $X_m = \xi_n$.
For this time $m$, the chance that $X_j \notin C_n$ for all $j > m$ is
equal to $\alpha$ by the strong Markov property.
Therefore, $\Psubbig_\bfid{\texists {m \ge n} \widehat
X_m \in \ewset_{n, \epsilon}(\Lamps_\infty)} \ge \alpha > 0$\rl{, as
desired}.

For the general case, let
$\tau_n$ be the first exit time of $C_n$ ($n \ge 0$).
Let $K_n$ be the ball of radius $R$ about $\xi_n$.
By transience, for each $x \in K_0$, there is some time $t_x \ge 0$ such
that $\Psub_x[\all {s \ge t_x} X_s \notin K_0] > 1/2$.
Choosing $\tmax := \max_{x \in K_0} t_x$ gives a time such that
$\Psub_x[\all {s \ge \tmax} X_s \notin K_0] > 1/2$ for all $x \in K_0$.
Before time $\tau_n$, a lamp can be changed only in $C_n \cup K_n$.
Let $A_n$ be the ball of radius $R (\tmax + 1)$ about $\xi_n$.
Then at times in $[\tau_n, \tau_n + \tmax]$,
the lamplighter must stay in $A_n$ and the
changes of lamps must be entirely within $A_n$.
We may define $\ewset_{n, \epsilon}(\lamps_\infty)$
to consist of those $(\lamps_n, x_n)$ such that $x_n \in A_{R n}$ and such that
$$
\lamps_n(y) =
\cases{\lamps_\infty(y) &if $y \in C_{R n} \setminus A_{R n}$,\cr
\id &if $y \notin C_{R n} \cup K_{R n} \cup A_{R n}$.\cr}
$$
Then $\ewset_{n, \epsilon}(\lamps_\infty)$ is of bounded size and
$\Psubbig_\bfid{\texists {m \ge n} \widehat X_m \in \ewset_{n,
\epsilon}(\Lamps_\infty)} \ge 1/2$.
\Qed

\rl{
Our last illustration of the enhanced criterion \ref c.vadim-quant2/
identifies the Poisson
boundary when the projection of $\mu$ on the
base group $\baseg$ does {\it not} generate
$\baseg$ as a {\it semigroup}.
Our proof in this case works for all nontrivial lamp groups and all
countably infinite base groups, not necessarily finitely generated.
We will not, however, need to use this result in our later proofs.
\procl t.not-semigp
Let $\lampg$ be a nontrivial group and $\baseg$ be an infinite group.
Let $\mu$ be a probability measure of finite entropy on $\lampg \wr \baseg$
whose support generates $\lampg \wr \baseg$ (as a group) and is
concentrated on $\big\{(\chngo {s}, o) \st s \in \lampg\big\} \cup
\big\{(\bfid, x) \st x \in \baseg\big\}$.
If the projection $\mub$ of $\mu$ on $\baseg$ 
has support that does {\rm not} generate $\baseg$ as a semigroup,
then the Poisson boundary
of $(\lampg \wr \baseg, \mu)$ is $\lampg^\baseg$ endowed with the law of\/
$\Lamps_\infty$.
\endprocl
The basic idea of the proof is that the random walk on the base group
has infinitely many cut times.
\proof
Let $\sgp$ denote the semigroup generated by the support of $\mub$,
including $o$.
Then $\sgp^{-1}$ is also a semigroup, as is
$\sgp' := \sgp \cap \sgp^{-1}$. Because $\sgp \cup \sgp^{-1}$
generates $\baseg$ as a semigroup,
$\sgp' \ne \sgp$ and 
$\alpha := \P[Y_1 \notin \sgp'] > 0$. Let $\bigl\{t \st
Y_t \notin \sgp' \bigr\}$ be listed as $\Seq{\tau_n \st n \ge 1}$ in
increasing order.
Note that for $x \in \sgp$ and $y \in \sgp \setminus \sgp^{-1}$, we have $x
y \sgp \subset \sgp \setminus \sgp^{-1}$.
Therefore, $X_{\tau_{n+1}} \sgp \subset X_{\tau_n} (\sgp \setminus \sgp^{-1})
\subsetneq X_{\tau_n} \sgp$. 
Furthermore, if $x \in (\sgp')^k$ for some $k \ge 0$, then 
$\sgp = x \sgp$ because $x \in \sgp'$.
That is, we have a monotonic decreasing sequence 
$$
\sgp
=
X_0 \sgp
=
X_1 \sgp
=
\cdots
=
X_{\tau_1 - 1} \sgp
\supsetneq
X_{\tau_1} \sgp
=
\cdots
=
X_{\tau_2 - 1} \sgp
\supsetneq
X_{\tau_2} \sgp
=
\cdots.
$$
\vskip\parskip
Given $x, y \in \sgp$, write $x \approx y$ if $x \sgp = y \sgp$, and write $x
\prec y$ if $x \sgp \supsetneq y \sgp$. Write $x \precsim y$ if
$x \approx y$ or $x \prec y$. Then for every $n$, we
have $s, t \in [\tau_n, \tau_{n+1})$ implies $X_s \approx X_t$, whereas if $s
< \tau_n \le t$, then $X_s \prec X_t$.
\vskip\parskip
Recall that $\lit \lamps$ denotes the set of lit lamps, $\{x \in \baseg \st
\lamps(x) \ne \id\}$, of $\lamps \in \lampg^\baseg$.
Define the stopping times $\sigma_n := \inf \bigl\{ t \st |\{s \le t \st
X_s \in \lit(\Lamps_\infty)\}| \ge n,\ X_t \in \lit(\Lamps_\infty),\
X_{t-1} \prec X_t \bigr\}$; necessarily, $\sigma_n \ge n$. Then
$\P[X_{\sigma_n} \prec X_{\sigma_n + 1}] \ge \alpha$.
On the event $[X_{\sigma_n} \prec X_{\sigma_n + 1}]$,
we have $\Lamps_{\sigma_n}(x) =
\Lamps_\infty(x)$ for all $x \in \{X_s \st s \le \sigma_n\}$ and
$\Lamps_{\sigma_n}(x) = \id$ for all other $x \in \gp$; also,
$X_t \approx X_{\sigma_n}$ only for $t = \sigma_n$ on that event.
\vskip\parskip
Let $\lamps_\infty \in \lampg^{\baseg}$ be a possible limiting lamp
configuration.
For every $x, y \in \lit(\lamps_\infty)$,
 exactly one of the following holds:
$x \approx y$, $\> x \prec y$, or $y \prec x$, because $x$ and $y$ lie in
the trace of the random walk on $\baseg$.
Define $Q_{n}(\lamps_\infty)$ to be
the set of all $(\lamps_n, x)$ such that
\beginitems
\item{(i)}
$x \in \lit(\lamps_\infty)$,
\item{(ii)}
$|\{y \precsim x \st y \in \lit(\lamps_\infty) \}| \ge n$,
\item{(iii)}
if $y \in \lit(\lamps_\infty)$ and $y \approx x$, then $y = x$,
\item{(iv)}
if $\lit(\lamps_\infty) \ni z \prec x$ and $|\{y \precsim z \st y \in
\lit(\lamps_\infty)\}| \ge n$,
then there is some $w \approx z$ with $w \ne z$ and $w \in \lit(\lamps_\infty)$,
\enditems
and
\beginitems
\item{(v)}
$\displaystyle
\lamps_n(y) =
\cases{\lamps_\infty(y)
  & if $y \precsim x$ and $y \in \lit(\lamps_\infty)$,
\vadjust{\kern2pt}\cr
 \id & otherwise.\cr}
$
\vadjust{\kern2pt}%
\enditems
We have shown that $\P[\texists {m \ge n}
\widehat X_m \in Q_{n}(\Lamps_\infty)] \ge \P[\widehat X_{\sigma_n} \in
Q_n(\Lamps_\infty)] \ge
\alpha$. In addition, $|Q_n(\Lamps_\infty)| \le 1$. Thus, the theorem
follows from \ref c.vadim-quant2/. 
\Qed
}

\bsection{Proof of Theorem 1.1}{s.jumps}

In this section, we prove \ref t.finiteent/. 
This comes in \rl{three} parts; one
handles base groups $\baseg$ that \rl{have at least cubic growth and}
are Liouville for the projected walk (\ref
t.finiteent-Liou/)\rl{; one handles base groups of less than cubic growth
(\ref t.finiteext/);} and
the \rl{last} handles the rest (\ref t.finiteent-nonLiou/).
\rl{In fact, \ref t.finiteent-Liou/ also handles some other cases; the
reader interested in those cases can thereby find a proof that is simpler
than the one that uses all three theorems.}
We will write ``with high probability" to mean ``with probability tending
to 1 as $n \to\infty$".

\rl{For ease in following our proofs, we will assume that $\mu$ is
concentrated on $\big\{(\chngo {s}, o) \st s \in \lampg\big\} \cup
\big\{(\bfid, x) \st x \in \baseg\big\}$.}
It will be easy to see that the same proofs---indeed, with
simplifications---extend to all $\mu$ whose
support is finite and generates $\lampg \wr \baseg$.
\rl{The extension to $\mu$ with bounded lamp range involves merely
technical complications.}

\rl{We begin with five short lemmas.}

\procl l.binomial
If $k \le n/3$, then $\sum_{j=0}^k {n \choose j} \le 2(ne/k)^k$.
\endprocl


\proof
Since $k! \ge (k/e)^k$ by Stirling's inequality (p.~54 of \ref
b.Feller:3rd/),
we have ${n \choose k} \le (ne/k)^k$. Since ${n \choose j+1} \ge 2 {n \choose
j}$ for $j < n/3$, the result follows by comparison with a geometric series.
\Qed

The following theorem of Shannon
is well known and easy to prove via the weak law of large
numbers (e.g., \ref b.CoverThomas/, Theorem 3.1.2).

\procl l.aep
If $\pd$ is a discrete distribution on a set $S$
with entropy $H(\pd)$ and $Y_n \sim \pd$ are
independent, then there are sets $\Lambda_n \subseteq S^n$ ($n \ge 1$) such
that $\lim_{n \to\infty} n^{-1} \log |\Lambda_n| = H(\pd)$ and $\lim_{n
\to\infty} \Pbig{(Y_1, Y_2, \ldots, Y_n) \in \Lambda_n} = 1$.
\Qed
\endprocl

Write $\dist(x, y)$ for the distance between $x$ and $y$ in some Cayley
graph of $\gp$ and
$V_\gp(r)$ for the number of points within distance $r$ of the identity,
$\oid$.
Let $B(x, r)$ denote the ball of radius $r$ about $x$.
\rl{The following lemma is well known in cases such as symmetric simple
random walk, due to celebrated results of Varopoulos; see,
e.g., Corollary 7.3 of \ref b.CGP/ or, for a short proof, Corollary 6.6 of
\ref b.LOS/. It is easily deduced for nonsymmetric random walks from known
results,
but for completeness, we include this derivation.}

\procl l.close-polyd
Let $\Seq{\walk_n}$ be a $\mu$-walk on a group $\gp$ that satisfies
$V_\gp(r) \ge c r^d$ for all $r \in \N$. Assume that the support of $\mu$
generates $\gp$ and that $\mu(\oid) \ge 1/2$. Then
$p_{t}(\oid, x) \le c t^{-d/2}$ for all $t \ge 1$ and all $x \in \gp$.
\endprocl

\proof
Let $P$ be the transition matrix for the $\mu$-walk and
$\widehat P$ denote its transpose, which is the transition matrix for
another random walk.
Since the support of $\mu$ generates $\gp$ and
$\mu(\oid) > 0$, some power $P^j$ has the property that $(P^j + \widehat
P^j)/2$ is irreducible.
Thus, for such $j$, we have $\alpha :=
\min \medbigl\{ p_j(x, y) + \widehat p_j(x, y) \st x \sim y \medbigr\} > 0$.
It is well known that $\sum_{x \in S,\ y \notin S} p_j(x, y) =
\sum_{x \in S,\ y \notin S} \widehat p_j(x, y)$
for all finite $S \subset \verts$ (e.g., see \ref b.MorrisPeres/).
Since the sum of these two quantities is at least
$\alpha |\{(x, y) \st x \in S,\ y \notin S\}|$, it follows that
$\sum_{x \in S,\ y \notin S} p_j(x, y) \ge \alpha |\{(x, y) \st x \in S,\ y
\notin S\}|/2$.
Now the result follows from the isoperimetric inequality of \ref b.CoulSC/
and Corollary 6.32(i) of \ref b.LP:book/.
\Qed

\vglue-\bigskipamount
\rl{
\procl l.summable
For every symmetric, transient $\mu$-walk on a group $\gp$, 
$$
\sum_{t=0}^\infty \sup_{x \in \gp} p_t(\oid, x) < \infty\,.
$$
\endprocl
\proof
It is well known that for even $t$, we have $p_t(\oid, x) \le p_t(\oid,
\oid)$. Choose $y$ with $\mu(y) > 0$.
For odd $t$, we have $p_t(\oid, x) \le p_{t+1}(\oid, xy)/p(x, xy) \le
p_{t+1}(\oid, \oid)/\mu(y)$. Thus, the result follows from
$\sum_{t=0}^\infty p_t(\oid, \oid) < \infty$. 
\Qed
}

\procl l.close-genl
Let $\Seq{\walk_n}$ be a random walk on $\gp$.
Let $0 \le k < m$.
Suppose that $\MM$ is a random subset of\/ $\gp$ that is measurable with
respect to $\Seq{\walk_1, \ldots, \walk_k}$.
Then
$$
\P[\walk_m \in \MM]
\le
\Ebig{|\MM|} \sup_{x \in \gp} p_{m - k}(\oid, x)
\,.
$$
\rl{If $\MM = B(\walk_k, r)$, then
$$
\P[\walk_m \in \MM]
=
\sum_{x \in B(\oid, r)} p_{m - k}(\oid, x)
\,.
$$}
\endprocl

\proof
For each $y \in \gp$, we have
$$
\P[\walk_m = y \mid \walk_1, \ldots, \walk_k]
=
p_{m - k}(\oid, \walk_k^{-1} y)
\le
\sup_{x \in \gp} p_{m - k}(\oid, x)
\,.
$$
Summing over $y \in \MM$ and then taking expectation gives the first result.
\rl{In the second case, we use instead the identity 
$$
\P[\walk_m \in \MM \mid \walk_1, \ldots, \walk_k]
=
\sum_{y \in \MM} p_{m - k}(\oid, \walk_k^{-1} y)
=
\sum_{x \in B(\oid, r)} p_{m - k}(\oid, x)
\,.
\Qed
$$}

\procl t.finiteent-Liou
Let $\lampg$ be a nontrivial finite group
\rl{and $\baseg$ be a finitely generated, infinite group.
Let $\mu$ be a probability measure of finite entropy on $\lampg \wr \baseg$
whose support generates $\lampg \wr \baseg$ and is
concentrated on $\big\{(\chngo {s}, o) \st s \in \lampg\big\} \cup
\big\{(\bfid, x) \st x \in \baseg\big\}$.
Suppose that the projection $\mub$ of $\mu$ on $\baseg$ is Liouville and
generates a transient random walk.
If any one of the following conditions holds,
then the Poisson boundary
of $(\lampg \wr \baseg, \mu)$ is $\lampg^\baseg$ endowed with the law of\/
$\Lamps_\infty$:
\beginitems
\itemrm{(a)} 
the measure $\mub$ is symmetric; or
\itemrm{(b)} 
the group $\baseg$ has at least cubic growth; or
\itemrm{(c)} 
the group $\baseg$ is abelian.
\enditems
\vskip0pt}
\endprocl

A rough sketch of the proof follows.
Since $\Seq{X_n}$ is Liouville, its asymptotic entropy is 0, whence there
is some $t_0$ such that $H(X_{t_0}) < \epsilon t_0$.
\ref l.aep/ converts this to a likely set of fewer than $e^{\epsilon n}$
possibilities for $S := \Seq{X_{j t_0} \st j \le n/t_0}$.
For a large $\rho$, partially obscure the increments $\Seq{Y_k \st 1 \le k
\le n}$ by replacing
those $Y_k$ that satisfy $|Y_k| \le \rho$ by $*$, to
mean ``unknown";
the resulting sequence $U \in \big(\baseg \cup \{*\}\big)^{n}$
has small entropy, so we again have a collection of size $< e^{\epsilon n}$
containing likely values $U$ of the partially obscured increments.
In this way, we guess the large jumps and bound the others.
Knowing $S$ and $U$, we define the set $M_i(S, U)$ of possible values for
$\{X_j \st (i-1)t_0 < j \le i t_0\}$, and $|M_i(S, U)| \le t_0 V_\baseg(\rho)^{t_0}$.
In most locations $y \in M_i(S, U)$, we have
$\Lamps_n(y) = \Lamps_\infty(y)$, and we
can bound the number of possibilities for $\Lamps_n(y)$ for the other $y$.

\proof
Since $H(X_1) < \infty$ and the walk on $\baseg$ is Liouville, we have
that $H(X_n) = o(n)$.

Let $\epsilon > 0$.
Choose $t_0$ so that $H(X_{t_0}) < \epsilon t_0$.
For $n \in t_0 \cdot \Z^+$, set $s_n := n/t_0$.
Write $\ranen S := \Seq{X_{j t_0} \st 1 \le j \le s_n}$.
Applying
\ref l.aep/ to the $t_0$-step increments $X_{j t_0}^{-1} X_{(j+1) t_0}$ yields
a set $\enum S_n \subseteq \baseg^{s_n}$ with $\log |\enum
S_n| < \epsilon t_0 s_n = \epsilon n$ and $\P[\ranen S \in \enum S_n] \to 1$.

Write
$$
u_\rho(x) := \cases{x &if\/ $\dist(o, x) > \rho$,\cr
              * &otherwise.\cr}
\label e.defu
$$
Recall that $\Seq{Y_k}$ are the increments of the random walk on $\baseg$.
Choose $\rho$ so that
$H\big(u_{\rho}(Y_1)\big) < \epsilon$.

Write $\ranen U := \Seq{u_\rho(Y_k) \st 1 \le k \le n}$.
By \ref l.aep/,
there is a set $\enum U_n \subseteq \big(\baseg \cup
\{*\}\big)^{n}$ with $\log |\enum U_n| < \epsilon n$ and
$\P[\ranen U \in \enum U_n] \to 1$.
For each $U \in \enum U_n$ and $0 \le j_1 < j_2 \le n$,
define the set $L(U, j_1, j_2) \subset \baseg$ to be the set of possible
values of $X_{j_1}^{-1} X_{j_2}$ that are consistent with $\ranen U = U$.
That is, let
$$
Z_k := \cases{
Y_k &if $|Y_k| > \rho$,\cr
B(o, \rho) &otherwise,\cr}
$$
and define
$$
L(U, j_1, j_2)
:=
\prod_{j_1 < k \le j_2} Z_k
:=
Z_{j_1+1} Z_{j_1+2} \cdots Z_{j_2}
\,.
\label e.defL
$$
When $\baseg$ is abelian, 
\rl{$L(U, j_1, j_2)$ is a ball of radius at most $\rho (j_2 - j_1)$.}
More generally,
$$
|L(U, j_1, j_2)|
\le
V_\baseg(\rho)^{j_2-j_1}
\le
V_\baseg(\rho)^{t_0}
\,.
$$
Given $S = \Seq{ x_1, x_2, \ldots, x_{s_n}}$ and $1 \le i \le s_n$,
write
$$
M_i(S, U)
:= \bigcup_{j=(i-1)t_0+1}^{i t_0} x_{i-1} L\bigl(U, (i-1) t_0, j\bigr)
\,,
$$
where $x_0 := o$, for the set of possible values of $\{X_j \st (i-1)t_0 < j
\le i t_0\}$ that are consistent with $\ranen S = S$ and $\ranen U = U$.
Thus, $\bigcup_{i \in [1, s_n]} M_i(S, U)$ \rl{contains all possible values of
$X_j$ for $0 < j \le n$}
that are consistent with $\ranen S = S$ and $\ranen U = U$; outside this
set, every lamp must be the identity at time $n$.
Inside this set, the lamp at time $n$ takes the same value as at time
$\infty$ except possibly at those locations that are visited after time
$n$.

Thus, let $\ranen W := \big\{i \in [1, s_n] \st \texists {m > n} X_m \in
M_i(\ranen S, \ranen U) \big\}$. At the end of the last paragraph, we
observed that
$$
\Lamps_n(z) =
\cases{\Lamps_\infty(z)
  & for $z \in \bigcup_{i \in [1, s_n] \setminus \ranen W} M_i(\ranen S, \ranen U)$,
\vadjust{\kern2pt}\cr
 \id & for $z \notin \bigcup_{i \in [1, s_n]} M_i(\ranen S, \ranen U)$.\cr}
$$
\rl{Choose $a_n := \sqrt{\Ebig{|\ranen W|}n}$. 
Write $A_n := \big[|\ranen W| \le a_n\big]$.
We claim that in all three
cases (a)--(c), $\lim_{n \to\infty} \Ebig{|\ranen W|}/n = 0$, whence
$\lim_{n \to\infty} \Ebig{|\ranen W|}/a_n = \lim_{n \to\infty} a_n/n = 0$.
Once we establish that, 
we may deduce that $\lim_{n \to\infty} \P(A_n) = 1$ by Markov's inequality.
\vskip\parskip
In order to show our claim, 
apply \ref l.close-genl/ to see that
for $m > n$, $\>1 \le i \le s_n$, and $(i-1)t_0 < j \le i t_0$,
$$\eqaln{
\Pbig{X_m \in X_{(i-1)t_0} L(\ranen U, (i-1)t_0, j)}
&\le
\Ebig{|L(\ranen U, (i-1)t_0, j)|} \sup_{x \in \baseg} p_{m-j}(o, x)
\cr &\le
V_\baseg(\rho)^{t_0} \sup_{x \in \baseg} p_{m-j}(o, x)
\,.}
$$
By virtue of Lemmas \briefref l.close-polyd/ and \briefref l.summable/, we
have that in cases (a) and (b), 
$$
\alpha_\ell := V_\baseg(\rho)^{t_0}
\sum_{k > \ell} \sup_{x \in \baseg} p_{k}(o, x)
\to 0
$$
as $\ell \to\infty$.
On the other hand, in case (c), $X_{(i-1)t_0} L(\ranen U, (i-1)t_0, j)$ is
a ball of radius at most $\rho t_0$ that contains $X_j$, whence
$X_{(i-1)t_0} L(\ranen U, (i-1)t_0, j) \subseteq B(X_j, 2\rho
t_0)$.
Thus, in case (c),
it follows from \ref l.close-genl/ that
for $m > n$, $\>1 \le i \le s_n$, and $(i-1)t_0 < j \le i t_0$,
$$
\Pbig{X_m \in X_{(i-1)t_0} L(\ranen U, (i-1)t_0, j)}
\le
\sum_{|x| \le 2\rho t_0} p_{m-j}(o, x)
\,.
$$
In this case, transience guarantees that
$$
\alpha_\ell := 
\sum_{k > \ell} \sum_{|x| \le 2\rho t_0} p_{k}(o, x)
\to 0
$$
as $\ell \to\infty$.
Therefore, in all three cases,
$$
\Ebig{|\ranen W|} 
\le
\sum_{j \le n}
\sum_{m> n}
\Pbig{X_m \in X_{(\ceil{j/t_0}-1)t_0} L(\ranen U, (\ceil{j/t_0}-1)t_0, j)}
\le
\sum_{j \le n} \alpha_{n - j}
=
\sum_{j \le n} \alpha_{j}
=
o(n)
$$
as $n \to\infty$, as claimed.
}

Let $\lamps_\infty \in \lampg^{\baseg}$ be a possible limiting lamp
configuration.
For $n \in t_0 \cdot \Z^+$, define $Q_{n, \epsilon}(\lamps_\infty)$ to be
the set of all
$(\lamps_n, x)$ such that
there are $U$, $S$, and $W$ satisfying
\beginitems
\item{(i)}
$\rl{U \in \enum U_n}$,
\item{(ii)}
$S = \Seq{ x_1, x_2, \ldots, x_{s_n}} \in \enum S_n$
with $x_{s_n} = x$,
\item{(iii)}
$W \subseteq [1, s_n]$ with $|W| \le \rl{a_n}$,
\enditems
and
\beginitems
\item{(iv)}
$\displaystyle
\lamps_n(z) =
\cases{\lamps_\infty(z)
  & for $z \in \bigcup_{i \in [1, s_n] \setminus W} M_i(S, U)$,
\vadjust{\kern2pt}\cr
 \id & for $z \notin \bigcup_{i \in [1, s_n]} M_i(S, U)$.\cr}
$
\vadjust{\kern2pt}%
\enditems

We have established that $\ranen U$,
$\ranen S$, and $\ranen W$ satisfy (i)--(iv) with high probability as
choices for $U$, $S$, and $W$, respectively, when
$\lamps_\infty = \Lamps_\infty$, $\>\lamps_n = \Lamps_n$, and $x = X_n$, and
thus $\lim_{t_0 \cdot \Z^+ \ni n \to\infty} \Pbig{\widehat
X_n \in Q_{n, \epsilon}(\Lamps_\infty)} = 1$. To establish the theorem, in
light of \ref c.Kcrit-epsilon/, it suffices to show that
$|Q_{n, \epsilon}(\lamps_\infty)| < e^{2 \epsilon n + o(n)}$ because
$\epsilon$ was arbitrary.

By definition, the number of choices of $U \in \enum U_n$ is at most $e^{
\epsilon n}$ and the number of choices of $S \in \enum S_n$ is at most $e^{
\epsilon n}$.
For \rl{large $n$},
the number of choices of $W$ is at most $2(\rl{{n \over a_n}e)^{a_n}} = e^{o(n)}$ by \ref
l.binomial/.
Note that $|M_i(S, U)| \le t_0 V_\baseg(\rho)^{t_0}$. Thus,
given $S$ and $W$, the number of choices of $\lamps_n$
is at most $|\lampg|^{|W| t_0 V_\baseg(\rho)^{t_0}} \le |\lampg|^{t_0
V_\baseg(\rho)^{t_0}\rl{a_n}} = e^{o(n)}$.
Therefore, $|Q_{n, \epsilon}(\lamps_\infty)| < e^{2 \epsilon n + o(n)}$, as
desired.
\Qed

\rl{
In order to handle the case of base groups with less than cubic growth, we
modify the preceding proof in a couple of ways. By \ref b.Gromov/,
all such groups are finite extensions of $\Z$ or $\Z^2$.
\procl t.finiteext
Let $\lampg$ be a nontrivial finite group and $\baseg$ be a finitely
generated, infinite group with an abelian
subgroup $\as$ of finite index.
Let $\mu$ be a probability measure of finite entropy on $\lampg \wr \baseg$
whose support generates $\lampg \wr \baseg$ and is
concentrated on $\big\{(\chngo {s}, o) \st s \in \lampg\big\} \cup
\big\{(\bfid, x) \st x \in \baseg\big\}$.
If the projection $\mub$ of $\mu$ on $\baseg$ 
generates a transient random walk,
then the Poisson boundary
of $(\lampg \wr \baseg, \mu)$ is $\lampg^\baseg$ endowed with the law of\/
$\Lamps_\infty$.
\endprocl
The idea of the proof is to use the commutativity of $\as$ to further
specify the
possible positions of the base walk in the first $n$ steps, beyond what the
previous proof accomplished. The aim is to pay most attention when the base
walk lies in $\as$. When the base walk moves far during an
excursion between visits to $\as$, then we will specify exactly the
increments during an entire such excursion. With ``far"
having a sufficiently large threshold, such specification can be done with
a collection of size $< e^{\epsilon n}$ of likely values.
There are extra difficulties because the times when the base walk lies in
$\as$ are random, but since we need to know only relatively few of them, we
can choose a possible set of such times with small exponential growth.
\proof
Note that if $\tau$ is a stopping time, then 
$$
H\bigl(\Seq{Y_t \st
1 \le t \le \tau}\bigr) = \Ebigg{\sum_{t=1}^\tau \log \mu(Y_t)} =
\E[\tau] \E[\log \mu(Y_1)] = \E[\tau] \, H(Y_1)
\,.
$$
In particular, this is finite when $\E[\tau] < \infty$.
\vskip\parskip
Since $H(X_1) < \infty$ and the $\mub$-walk on $\baseg$ is necessarily
Liouville, we have that $H(X_n) = o(n)$.
Let $\epsilon > 0$.
Choose $t_0$ so that $H(X_{t_0}) < \epsilon t_0$.
For $n \in t_0 \cdot \Z^+$, set $s_n := n/t_0$.
Write $\ranen S := \Seq{X_{j t_0} \st 1 \le j \le s_n}$.
Applying
\ref l.aep/ to the $t_0$-step increments $X_{j t_0}^{-1} X_{(j+1) t_0}$ yields
a set $\enum S_n \subseteq \baseg^{s_n}$ with $\log |\enum
S_n| < \epsilon t_0 s_n = \epsilon n$ and $\P[\ranen S \in \enum S_n] \to 1$.
\vskip\parskip
For a sequence $\overline x = \Seq{x_j \st 1 \le j \le t}$, write
$$
\dmax(\overline x) := \max \{
\dist(o, x_j) \st 1 \le j \le t\}\,.
$$
Write 
$$
u_\rho(\overline x)
:=
\cases{\overline x &if\/ $\dmax(\overline x) > \rho$,\cr
       * &otherwise.\cr}
$$
Let $\tau_\as(k)$ be the time of the $k$th visit of the $\mub$-walk 
to $\as$, with $\tau_\as(0) := 0$.
Because $[\baseg : \as]$ is finite, $\Ebig{\tau_\as(1)} < \infty$.
Abbreviate the sequence $\Seq{X_{\tau_\as(k-1)}^{-1} X_j \st
\tau_\as(k-1) < j \le \tau_\as(k)}$
as $\overline X_k$.
Thus, $\overline X_k$ are IID with finite entropy by our first paragraph.
Choose $\rho$ so that $H\bigl(u_{\rho}(\overline X_1)\bigr) < \epsilon$.
We may also assume that $\dist(x, \as) \le \rho$ for all $x \in \baseg$.
\vskip\parskip
Write 
$
\ranen U := \Seqbig{u_\rho(\overline X_k) \st 1 \le k \le n}
$.
Write $\baseg^{<\infty}$ for the set of finite sequences of
elements from $\baseg$.
By \ref l.aep/,
there is a set $\enum U_n \subseteq \big(\baseg^{<\infty}\cup
\{*\}\big)^{n}$ with $\log |\enum U_n| < \epsilon n$ and
$\P[\ranen U \in \enum U_n] \to 1$.
\vskip\parskip
The times $\tau_\as(k)$ for which $u_\rho(\overline X_k) \ne *$ form a
renewal process.
Let $\ranen T$ be the set of such renewal times $\le n$.
The long-term rate $\alpha_\rho$ of renewals tends to 0 as $\rho
\to\infty$.
Let $\enum T_n$ be the collection of subsets of $\{0, 1, 2,
\dots, n\}$ with size at most $2\alpha_\rho n$.
We have $\P[\ranen T \in \enum T_n] \to 1$ as $n \to\infty$.
By \ref l.binomial/, for sufficiently large $\rho$, the size of $\enum T_n$
is less than $e^{\epsilon n}$ for all large $n$.
Without loss of generality, we may assume that
$\rho$ is that large.
Let also $\ranen T'$ be the set of times $\tau_\as(k-1) \le n$ for which 
$u_\rho(\overline X_k) \ne *$.
Because $|\ranen T'| \le|\ranen T| + 2$,
we also have $\P[\ranen T' \in \enum T_n] \to 1$ as $n \to\infty$.
Write $\enum T''_n$ for the set of pairs $(T, T') \in \enum T_n \times
\enum T_n$ that are possible values of $(\ranen T, \ranen T')$: that is,
they must interleave with the minimum coming from $T'$.
\vskip\parskip
Observe that $(\ranen S, \ranen T, \ranen T', \ranen U)$ determines
$X_j$ for $j \in [1, s_n] t_0$ and also for $(i-1) t_0 \le j \le \tau_\as(k)$
when $\tau_\as(k-1) \le (i-1) t_0$ and $u_\rho(\overline X_k) \ne *$, where
$1 \le i \le s_n$. 
For all other $j \le n$, that quadruple forces $X_j$ to lie in a ball of radius
$3\rho t_0$ about some point $\beta_j \in \as$
that is measurable with respect to
$(\ranen S, \ranen T, \ranen T', \ranen U)$.
Indeed, fix a map $\zeta \colon \baseg \to \as$ such that $\dist\bigl(x,
\zeta(x)\bigr) \le \rho$.
Let $j \in \bigl[(i-1) t_0, i t_0\bigr)$.
Define $k$ by
$\tau_\as(k-1) < (i-1) t_0 \le \tau_\as(k)$.
Suppose first that $j \le \tau_\as(k)$.
In case $u_\rho(\overline X_k) \ne *$, then $X_j$ is determined by
$(\ranen S, \ranen T, \ranen T', \ranen U)$, so we may take $\beta_j :=
\zeta(X_j)$, whereas if
$u_\rho(\overline X_k) = *$, then $X_j \in B\bigl(\zeta(X_{(i-1) t_0}),
2\rho\bigr)$, so we may take $\beta_j := \zeta(X_{(i-1)t_0})$.
Suppose next that
$j > \tau_\as(k)$.
Write $\xi$ for the product of $X^{-1}_{\tau_\as(\ell-1)}
X_{\tau_\as(\ell)}$ over all $\ell > k$ with $\tau_\as(\ell) \le j$ and
$u_\rho(\overline X_\ell) \ne *$.
Because $\as$ is abelian, $X_j \in B\bigl(X_{\tau_\as(k)} \xi, \rho t_0\bigr)$
and $X_{\tau_\as(k)} \xi$ is $(\ranen S, \ranen T, \ranen T', \ranen
U)$-measurable
in case $u_\rho(\overline X_k) \ne *$, so we may take $\beta_j :=
\zeta(X_{\tau_\as(k)})$, and
$X_j \in B\bigl(\zeta(X_{(i-1) t_0}) \xi,
2\rho + \rho t_0\bigr)$
and $X_{(i-1) t_0} \xi$ is $(\ranen S, \ranen T, \ranen T', \ranen
U)$-measurable in the other case, so we may take $\beta_j :=
\zeta(X_{(i-1)t_0}) \xi$.
\vskip\parskip
For each $S \in \enum S_n$, $\>(T, T') \in \enum T''_n$, $\>U \in \enum
U_n$, $\>i \in [1, s_n]$, and $j \in \bigl((i-1)t_0, i t_0\bigr]$, define
the set $L(S, T, T', U, j) \subset \baseg$ to be the set of possible values
of $X_{j}$ that are consistent with $\ranen S = S$,
$\>\ranen T = T$, $\>\ranen T' = T'$, and $\ranen U = U$.
The preceding paragraph established that
$L(S, T, T', U, j)$ is contained in a ball of radius $3 \rho t_0$.
Write
$$
M_i(S, T, T', U)
:= \bigcup_{j=(i-1)t_0+1}^{i t_0} L(S, T, T', U, j)
\,.
$$
Thus, $\bigcup_{i \in [1, s_n]} M_i(S, T, T', U)$ contains all possible
values of $X_j$ for $0 < j \le n$
that are consistent with $(\ranen S, \ranen T, \ranen T', \ranen U) = (S,
T, T', U)$; outside this
set, every lamp must be the identity at time $n$.
Inside this set, the lamp at time $n$ takes the same value as at time
$\infty$ except possibly at those locations that are visited after time
$n$.
\vskip\parskip
Thus, let $\ranen W := \big\{i \in [1, s_n] \st \texists {m > n} X_m \in
M_i(\ranen S, \ranen T, \ranen T', \ranen U) \big\}$. At the end of the last paragraph, we
observed that
$$
\Lamps_n(z) =
\cases{\Lamps_\infty(z)
  & for $z \in \bigcup_{i \in [1, s_n] \setminus \ranen W} M_i(\ranen S,
  \ranen T, \ranen T', \ranen U)$,
\vadjust{\kern2pt}\cr
 \id & for $z \notin \bigcup_{i \in [1, s_n]} M_i(\ranen S, \ranen T,
 \ranen T', \ranen U)$.\cr}
$$
Choose $a_n := \sqrt{\Ebig{|\ranen W|}n}$. 
Write $A_n := \big[|\ranen W| \le a_n\big]$.
We claim that 
$\lim_{n \to\infty} \Ebig{|\ranen W|}/n = 0$, whence
$\lim_{n \to\infty} \Ebig{|\ranen W|}/a_n = \lim_{n \to\infty} a_n/n = 0$.
Once we establish that, 
we may deduce that $\lim_{n \to\infty} \P(A_n) = 1$ by Markov's inequality.
\vskip\parskip
In order to show our claim, 
apply \ref l.close-genl/ to see that
for $m > n$, $\>1 \le i \le s_n$, and $(i-1)t_0 < j \le i t_0$,
because 
$L(\ranen S, \ranen T, \ranen T', \ranen U, j) \subseteq B\bigl(X_j, 6\rho 
t_0\bigr)$,
$$
\Pbig{X_m \in L(\ranen S, \ranen T, \ranen T', \ranen U, j)}
\le
\sum_{|x| \le 6\rho t_0} p_{m-j}(o, x)
\,.
$$
Because of transience,
$$
\alpha_\ell := 
\sum_{k > \ell} \sum_{|x| \le 6\rho t_0} p_{k}(o, x)
\to 0
$$
as $\ell \to\infty$,
whence 
$$
\Ebig{|\ranen W|} 
\le
\sum_{j \le n}
\sum_{m> n}
\Pbig{X_m \in L(\ranen S, \ranen T, \ranen T', \ranen U, j)}
\le
\sum_{j \le n} \alpha_{n - j}
=
\sum_{j \le n} \alpha_{j}
=
o(n)
$$
as $n \to\infty$, as claimed.
\vskip\parskip
Let $\lamps_\infty \in \lampg^{\baseg}$ be a possible limiting lamp
configuration.
For $n \in t_0 \cdot \Z^+$, define $Q_{n, \epsilon}(\lamps_\infty)$ to be
the set of all
$(\lamps_n, x)$ such that
there are $S$, $T$, $T'$, $U$, and $W$ satisfying
\beginitems
\item{(i)}
$S = \Seq{ x_1, x_2, \ldots, x_{s_n}} \in \enum S_n$
with $x_{s_n} = x$,
\item{(ii)}
$(T, T') \in \enum T_n''$,
\item{(iii)}
$U \in \enum U_n$,
\item{(iv)}
$W \subseteq [1, s_n]$ with $|W| \le \rl{a_n}$,
\enditems
and
\beginitems
\item{(v)}
$\displaystyle
\lamps_n(z) =
\cases{\lamps_\infty(z)
  & for $z \in \bigcup_{i \in [1, s_n] \setminus W} M_i(S, T, T', U)$,
\vadjust{\kern2pt}\cr
 \id & for $z \notin \bigcup_{i \in [1, s_n]} M_i(S, T, T', U)$.\cr}
$
\vadjust{\kern2pt}%
\enditems
\vskip\parskip
We have established that $(\ranen S, \ranen T, \ranen T', \ranen U, \ranen
W)$ satisfy (i)--(v) with high probability as a 
choice for $(S, T, T', U, W)$ when
$\lamps_\infty = \Lamps_\infty$, $\>\lamps_n = \Lamps_n$, and $x = X_n$, and
thus 
$$
\lim_{t_0 \cdot \Z^+ \ni n \to\infty} \Pbig{\widehat
X_n \in Q_{n, \epsilon}(\Lamps_\infty)} = 1
\,.
$$
To establish the theorem via
\ref c.Kcrit-epsilon/, it suffices to show that
$|Q_{n, \epsilon}(\lamps_\infty)| < e^{4 \epsilon n + o(n)}$ since
$\epsilon$ was arbitrary.
\vskip\parskip
By definition, 
$|\enum S_n| < e^{ \epsilon n}$, 
$\>|\enum T_n''| < e^{ 2\epsilon n}$, and
$|\enum U_n| < e^{ \epsilon n}$.
For large $n$,
the number of choices of $W$ is at most $2({n \over a_n}e)^{a_n} =
e^{o(n)}$ by \ref l.binomial/.
Note that $|M_i(S, T, T', U)| \le t_0 V_\baseg(3\rho t_0)$. Thus,
given $(S, T, T', U, W)$, the number of choices of $\lamps_n$
is at most 
$$
|\lampg|^{|W| t_0 V_\baseg(3\rho t_0)} \le |\lampg|^{t_0
V_\baseg(3\rho t_0)a_n} = e^{o(n)}\,.
$$
Therefore, $|Q_{n, \epsilon}(\lamps_\infty)| < e^{4 \epsilon n + o(n)}$, as
desired.
\Qed
}

\procl t.finiteent-nonLiou
Let $\lampg$ be a nontrivial \rl{finite or countable} group and
$\baseg$ be a \rl{countably infinite} group.
Let $\mu$ be a probability measure of finite entropy on $\lampg \wr \baseg$
whose support generates $\lampg \wr \baseg$ and that is
concentrated on $\big\{(\chngo {s}, o) \st s \in \lampg\big\} \cup
\big\{(\bfid, x) \st x \in \baseg\big\}$.
If the projection of $\mu$ on $\baseg$ is non-Liouville, then the Poisson boundary
of $(\lampg \wr \baseg, \mu)$ is $\lampg^\baseg$ endowed with the law of\/
$\Lamps_\infty$.
\endprocl

A rough sketch of the proof follows.
Let $\tau_x := \inf \{ n \st X_n = x\}$.
Let $\grnmet_\baseg(x) := -\log \Psub_o[\tau_x < \infty]$, the negative log of
the probability that the projection of the $\mu$-walk to $\baseg$, started at
$o$, ever visits $x \in \baseg$.
Because the walk on $\baseg$ is non-Liouville, its Avez entropy is $\ent' >
0$. It is known that $n^{-1} \grnmet_\baseg(X_n) \to \ent'$ a.s.
Consider the sets $W(r) := \{x \in \baseg \st \grnmet_\baseg(x) \le r\}$,
the sizes of which will not concern us.
Given $\epsilon > 0$, it is likely that for large $n$, we have $X_k \in
W := W\big(n \ent'(1 + \epsilon)\big)$ for all $k \le n$ and also that $X_m
\notin W$ for all $m > n(1+3\epsilon)$.
At the same time, there is a reasonable chance that $\Lamps_\infty(X_n) \ne
\id$.  Thus, there is a reasonable chance that $\Lamps_n$ agrees with
$\Lamps_\infty$ on $W \setminus \{X_{n+1}, \ldots, X_{n(1+\epsilon)}\}$,
and it is likely that
$\Lamps_n(z) = \id$ for all $z \notin W$. Furthermore, there are likely fewer
than $n(1+3\epsilon)$ locations $z \in W$ where $\Lamps_\infty(z)\ne \id$,
whence it is likely, seeing $\Lamps_\infty$,
that there are not many possibilities
for where $X_n$ is. Finally, it is likely that
$\widehat Y_m$ for $n < m \le
n(1+3\epsilon)$ belongs to a set of size $e^{c \epsilon n}$ (that does not
depend on $\Lamps_\infty$). From $\Lamps_\infty$ and these possibilities,
we can thus likely deduce
$\Lamps_n$.

\proof
Let $\ent' > 0$ be the Avez entropy of the projection of the $\mu$-walk to
$\baseg$.
By Proposition 6.2 of
\ref b.BP:tirw/ in the symmetric case
or \ref b.BHM/ in general, $\lim_{n \to\infty} n^{-1}
\grnmet_\baseg(X_n) = \ent'$ a.s.; this result is also proved as Theorem 14.50
of \ref b.LP:book/.
Write $W(r) := \{x \in \baseg \st \grnmet_\baseg(x) \le r\}$.
Let $\epsilon \in (0, 1/3)$.
Let $W := W\big(n \ent'(1 + \epsilon)\big)$
and $W' := W\big(n(1 + 3\epsilon) \ent' (1 - \epsilon)\big)$.
Since $(1+3\epsilon)(1-\epsilon) - (1+\epsilon) = \epsilon (1 - 3\epsilon) >
0$, it follows that $W' \supset W$.

Write $\ranen U := \Seq{\widehat Y_m \st n < m \le n + 3 \epsilon n}$.
By \ref l.aep/,
there is a set $\enum U_n \subseteq \big(\lampg \wr \baseg \big)^{\flr{3
\epsilon n}}$
with $\log |\enum U_n| < 6 \epsilon n H(\widehat X_1)$ and
$\P[\ranen U \in \enum U_n] \to 1$ as $n \to\infty$.

Let $A_n$ be the event that $X_k \in W$ for all $k \le n$.
Since every sequence $\Seq{r_n}$ with $\lim_{n \to\infty} r_n/n = \ent'$ has
the property that
for all sufficiently large $n$,
$\all {k \le n} r_k < n \ent' (1+\epsilon)$,
we have $\lim_{n \to\infty} \P(A_n) = 1$.
In addition, at any time,
the walk may leave its current location with the lamp not equal to $\id$,
after 1 or 2 steps, and never return.
Therefore, $\inf_n \P[\Lamps_\infty(X_n) \ne \id] > 0$.

Recall that $\lit \lamps$ denotes the set of lit lamps, $\{x \in \baseg \st
\lamps(x) \ne \id\}$, of $\lamps \in \lampg^\baseg$.
Let $D_n$ be the event that $X_m \notin W$ for all $m > n + 3 \epsilon
n$.
Since $\Pbig{\all{m > n} X_{m} \notin W\big(n(1+3\epsilon) \ent' (1 -
\epsilon)\big)} \to 1$ as $n \to\infty$ and $W\big(n(1+3\epsilon) \ent' (1
- \epsilon)\big) = W' \supset W$, it follows that $\P(D_n) \to 1$.
On the event $D_n$, we have that $|W \cap \lit \Lamps_\infty| \le n +
3 \epsilon n$.
Now, the lamp at any $z \in W$ is changed at time $m \in (n,
n+3\epsilon n]$ by multiplying by $\lincr_m(X_{m-1}^{-1} z)$, whence the
total change from what it was at time $n$ due to the changes in $\ranen U$ is
$\prod_{m=n+1}^{n + \flr{3 \epsilon n}} \lincr_m(X_{m-1}^{-1} z)$.
Therefore, on the event $D_n$, we have that for every $z \in W$,
$$
\Lamps_\infty(z)
=
\Lamps_n(z) \prod_{m=n+1}^{n + \flr{3 \epsilon n}} \lincr_m(X_{m-1}^{-1} z)
\,.
$$

%
Let $\lamps_\infty \in \lampg^{\baseg}$.
Define $Q_{n, \epsilon}(\lamps_\infty)$ to be the set of all
$(\lamps_n, x)$ such that
there is $U$ satisfying
\beginitems
\item{(i)}
$U = \Seq{(\psi_{n+1}, y_{n+1}), \ldots, (\psi_{n+\flr{3 \epsilon n}}, y_{n
+ \flr{3 \epsilon n}})}
\in \enum U_n$,
\item{(ii)}
$|W \cap \lit \lamps_\infty| \le n + 3 \epsilon n$,
\item{(iii)}
$x \in W \cap \lit \lamps_\infty$,
\enditems
and
\beginitems
\item{(iv)}
writing $z_m := x \prod_{j=n+1}^m y_j$ for $n \le m \le n + \flr{3 \epsilon
n}$ and
$$
\psi(z) :=
\prod_{m=n+1}^{n + \flr{3 \epsilon n}} \psi_m(z_{m-1}^{-1} z)
\,,
$$
we have
$$
\lamps_n(z) = \cases{\lamps_\infty(z)\psi(z)^{-1} & for $z \in W$,\cr
               \id & for $z \notin W$.\cr}
$$
\enditems

Clearly using $U := \ranen U$ satisfies (i) with high probability.
We have proved that for $\lamps_\infty = \Lamps_\infty$ and $x = X_n$, the
probability of (iii) is bounded away from 0. In addition,
(ii) and (iv) hold on the event $D_n$, which is likely.
Therefore, $\limsup_{n \to\infty} \Pbig{\widehat X_n \in
Q_{n, \epsilon}(\Lamps_\infty)} > 0$.

By assumption, the number of choices of $U \in \enum U_n$ is at most $e^{
6 \epsilon n H(\widehat X_1)}$.
The number of choices of $x$ is at most $2 n$.
Therefore, $|Q_{n, \epsilon}(\lamps_\infty)| < e^{6 \epsilon n H(\widehat
H_1) + o(n)}$.
This completes the proof.
\Qed

%

Define $\zeta_n(x) := -\log \Psub_o[\tau_x \le n]$.
We remark that one may use in the proof the more elementary fact that
$$
\lim_{n \to\infty} - n^{-1} \zeta_n(X_n) = \ent'
$$
(\ref b.BP:tirw/, proof of Proposition 6.2) in place of
$\lim_{n \to\infty} n^{-1} \grnmet_\baseg(X_n) = \ent'$.

\rl{
Lastly, we explain why a recurrent base walk yields a Liouville measure.
A group is called \dfn{Choquet--Deny} if every convolution walk on it is
Liouville. \ref b.FHTVF/ prove that such groups are exactly those groups
with no ICC quotients, where an \dfn{ICC} group is a nontrivial group all
of whose elements other than the identity have infinite conjugacy classes.
\vskip\parskip
\procl p.recurrent
Let $\lampg$ be a Choquet--Deny group and
$\baseg$ be a countable group.
Let $\mu$ be a probability measure on $\lampg \wr \baseg$
whose support generates $\lampg \wr \baseg$.
If $\mub$
generates a recurrent random walk on $\baseg$, 
then $(\lampg \wr \baseg, \mu)$ is Liouville.
\endprocl
\vskip\parskip
\proof
By the assumption that the $\mub$-walk is recurrent, the subgroup $\Delta$ of
elements of the form $(\Lamps, o)$ is a recurrence set. Let $\nu$ denote
the probability measure giving the first return to $\Delta$ from the identity.
Then the Poisson boundary of $(\lampg \wr \baseg, \mu)$ is isomorphic to
that of $(\Delta, \nu)$ by Lemma 4.2 of \ref b.Furst:Lie/.
Clearly, $\Delta$ is isomorphic to a direct sum of copies of $\lampg$.
On the other hand,
the direct sum of Choquet--Deny groups is Choquet--Deny. To see this,
let $\gp_i$ ($i \ge 1$) each have no ICC quotients and $\gp$ be their
direct sum.
Identify $\gp_i$ with the subgroup of elements of $\gp$ all of whose
coordinates are the identity other than the $i$th.  Let $\phi$ be a
homomorphism of $\gp$. We want to show that $\phi(\gp)$ is not ICC. Because
$\gp$ is generated by all $\gp_i$, we have that $\phi(\gp)$ is generated by
all $\phi(\gp_i)$. If $\phi(\gp)$ is not trivial, then some $\phi(\gp_i)$ is
not trivial and, by hypothesis, has a nontrivial element $\phi(\gpe_i)$, where
$\gpe_i \in \gp_i$, with finite conjugacy class in $\phi(\gp_i)$.  Since
$\gp_i$ commutes with all other $\gp_j$ ($j \ne i$) and $\phi$ is a
homomorphism, the conjugacy class of $\phi(\gpe_i)$ in $\phi(\gp_i)$ is the
same as in $\phi(\gp)$. Thus, we obtain our desired result that $\phi(\gp)$
contains a
nontrivial element with finite conjugacy class, so is not ICC. 
\Qed
}

\bsection{General Generators}{s.gen}

Here we extend the result of \ref b.Erschler/ from finite third moments to
finite second moments on $\Z^d$, and from $d \ge 5$ to $d \ge 3$.
We also allow infinite lamp groups.

\procl t.2ndmoment
Let $\lampg$ be a nontrivial finitely generated group and $d \ge 3$.
Let $\mu$ be a probability measure on $\lampg \wr \Z^d$
whose support generates $\lampg \wr \Z^d$
with
$\sum_{x} |x|^2 \,\mu(x) < \infty$.
Then the Poisson boundary of
$(\lampg \wr \Z^d, \mu)$ is $\lampg^\baseg$ endowed with the law of\/
$\Lamps_\infty$.
\endprocl

Note that for $x \in \Z^d$, its graph distance $|x|$ to $\bfz \in \Z^d$
is comparable to the $\ell^2$-norm $\|x\| :=
\iprod{\Cov(X_1) x}{x}^{1/2}$, which we define for $x \in \R^d$.
Write $B(r) := \{z \st \|z\| \le r\}$.

\rl{We have assumed that $\lampg$ is finitely generated only for brevity in
the assumptions; see the first paragraph of the proof for what we use
without this assumption.}

We preface the proof of \ref t.2ndmoment/ with a sketch.
The case when $\E[X_1] \ne \bfz$ was established by \ref b.Kaim:survey/, so
assume that $\E[X_1] = \bfz$.
The main new difficulty compared to our previous proofs
is that lamps may be changed at distances
arbitrarily far from the lamplighter. Control over this distance is
given by the moment assumption.
When $s$ is a large constant, for each $n$
there is a high chance that the first $n$ steps of
the walk on the base $\Z^d$ do not
exit the ball $B(s\sqrt n\,)$, nor
change any lamps outside the ball $B(2s \sqrt n\,)$.
In particular, there are only $c n^{d/2}$ possibilities for $X_n$ in this
case.
There is a tiny, but bounded below, chance that the walk on $\Z^d$ also has the
property that it never visits
the ball $B(4s \sqrt n\,)$ after time $n(1+\epsilon)$; conditional on this
event, the chance is very small that any lamp in $B(2s
\sqrt n\,)$ is changed after time $n(1+\epsilon)$.
There is a set of size $e^{c \epsilon n}$ that is likely to
contain $\widehat Y_{n+1}, \ldots, \widehat Y_{n(1+\epsilon)}$.
Having guessed $X_n \in B(s \sqrt n\,)$, seeing $\Lamps_\infty \restrict B(s
\sqrt n\,)$, and having changed the lamps therein according to
$\widehat Y_{n+1}, \ldots, \widehat Y_{n(1+\epsilon)}$, we arrive at
our guess of $\widehat X_n$.

\procl l.polylamps
Let $d \ge 3$.
Consider a random walk $\Seq{X_n}$ on $\Z^d$ with
$\Ebig{|X_1|^2} < \infty$ and $\E[X_1] = \bfz$.
\beginitems
\itemrm{(i)}
We have $\lim_{s \to\infty} \inf_n \P_\bfz \big[\all {k \le n} \|X_k\| \le s
\sqrt n\,\big] = 1$.
\itemrm{(ii)}
For every $s > 0$,
$$
\lim_{n \to\infty} \inf_{\|x\| \ge 2s\sqrt n} \Psubbig_x{\all {m \ge 0}
\|X_m\| > s \sqrt n \,}
= 1 - {1 \over 2^{d-2}}
\,.
$$
\enditems
\endprocl

\proof
Part (i) is immediate from Kolmogorov's maximal inequality (Theorem 2.5.2 in
\ref b.Durrett:book/).

To prove part (ii), we let $A:=\Cov(X_1)^{-1/2}$ and define
$Y_m:=s^{-1}AX_m$. Let $|\cdot|_2$ be the standard Euclidean norm and
$B_2(r)$ be the associated closed ball of radius $r$ about the origin.
Then (ii) can be rewritten in the form
$$
\lim_{n \to\infty} \sup_{|y|_2 \ge 2\sqrt n} \Psubbig_y{\texists  {m \ge 0}
|Y_m|_2 \le   \sqrt n \,} =2^{d-2}\,.
$$

First recall that
if standard Brownian in $\R^d$ starts at $z$ with $|z|_2=2$,
then the probability that it ever visits the ball $B_2(\alpha)$ of
radius $\alpha<2$ is $(\alpha/2)^{d-2}$; see, e.g., \ref b.BMbook/,
Corollary 3.19.
Given $\epsilon>0$, we can select $T=T(\epsilon)$  so that the probability
this visit occurs before time $T$ is at least $(\alpha/2)^{d-2}-\epsilon$;
taking $\alpha:=1-\epsilon$, we deduce from
the $d$-dimensional Donsker invariance
principle (see, e.g., \ref b.Whitt/,
Theorem 4.3.5),
$$
\lim_{n \to\infty} \sup_{|y|_2 \ge 2\sqrt n} \Psubbig_y{\texists  {m \in [0,nT]}
|Y_m|_2 \le  \sqrt n\, }
\ge  \bigl((1 - \epsilon)/2\bigr)^{d-2}-\epsilon
\,,
$$
and this gives the lower bound in (ii).

In dimension three, a matching upper bound follows from the asymptotic relation
$$
g(0,y) =\bigl(c_3+o(1)\bigr) |y|_2^{-1}
$$
for the Green function of $\Seq{Y_j}$ (see \ref b.Spitzer/, Proposition
P26.1), where $c_3$ is a positive constant. Indeed, if $\tau$ is the
hitting time of the ball $B_2(\sqrt{n}\,)$ by $\Seq{Y_j}$ (which may be
infinite), then the
optional stopping theorem (e.g., \ref b.Durrett:book/, Theorem 5.7.4)
for the bounded martingale $\Seq{g(0,Y_{\tau \wedge j}) \st j \ge 0}$ yields
$$
g(0,y)\ge \Psub_y [\tau<\infty] \cdot \min_{y_1 \in B_2(\sqrt{n}\,)} g(0,y_1)
\,.
$$
It follows that
$$
\bigl(c_3+o(1)\bigr) |y|_2^{-1}  \ge \Psub_y [\tau<\infty]
\bigl(c_3+o(1)\bigr) n^{-1/2}\,;
$$
the two occurrences of $o(1)$ do not necessarily denote the same function.
Since $|y|_2 \ge 2\sqrt n$, we conclude that
$\Psub_y [\tau<\infty]  \le 1/2+o(1)$.

It remains to prove the upper bound in (ii) for dimensions $d>3$. Given
$\epsilon>0$, let $T:=\ceil{d \epsilon^{-4}}$.
Another application of Donsker's theorem yields
$$
\lim_{n \to\infty} \sup_{|y|_2 \ge 2\sqrt n} \Psubbig_y{\texists  {m \in [0,nT]}
|Y_m|_2 \le   \sqrt n \,}
\le  (1/2+\epsilon)^{d-2}
\,.
$$
By the central limit theorem, for every $y$ in $\R^d$ and sufficiently large $n>1$,
$$
\Psubbig_y{|Y_{nT}|_2 \le \epsilon \sqrt{nT}\,} \le C_d \epsilon^d \,.
$$
If $z\in \R^d$ satisfies $|z|_2 > \epsilon \sqrt{nT} \ge \sqrt{nd}/\epsilon$,
then one of the coordinate projections of $z$ must have absolute value
greater than $\sqrt{n}/\epsilon$; projecting to a three-dimensional space
containing that coordinate, we infer (from optional stopping in three
dimensions) that as $ n \to \infty$,
$$
\Psubbig_z {\texists  {m \in [0,\infty)}  |Y_m|_2 \le   \sqrt n \,} \le
\bigl(1+o(1)\bigr)\epsilon \,.
$$
By considering whether $|Y_{nT}|_2 \le \sqrt{nd}/\epsilon$, we conclude that
for some constant $C_d$,
$$
 \sup_{|y|_2 \ge 2\sqrt n} \Psubbig_y{\texists  {m \ge 0} |Y_m|_2 \le
 \sqrt n\, } \le    (1/2+\epsilon)^{d-2}+C_d\epsilon^d + \bigl(1+o(1)\bigr)\epsilon \,.
\Qed
$$

Let $\rad$ denote the radius of a subset of $\baseg$, meaning the maximum
distance in the word metric of any of its elements from the identity
$o \in \baseg$.
Although we will apply the following lemma only for $\baseg = \Z^d$, we state
it in general as it may find other uses.
This lemma controls the changes of lamps far from the projection of the
walk on the base group.

\procl l.separate
Let $\lampg$ \rl{be a group and $\baseg$ be a finitely generated group.}
Let $\Seq{(\lincr_k, Y_k) \st k \ge 1}$ be the increments of a $\mu$-walk
$\Seq{\widehat X_n \st n \ge 0}$ on $\lampg \wr \baseg$
such that $\P[Y_1 = o] \ge 1/2$.
Suppose that $V_\baseg(r)/r^d$ is bounded above and below by positive finite constants
for some $d \ge 3$
and that $\Ebig{(\rad \lit \lincr_1)^2} < \infty$.
Then for some constant $c_\mu$, we have for every $a > 0$ that
$$
\sum_{k \ge 1} \Pbig{\rad \lit \lincr_k > a \cdot |X_{k-1}|}
\le
c_\mu a^{-2} \Ebig{(\rad \lit \lincr_1)^2}
\,.
$$
\endprocl

\proof
Let $R$ be a random variable independent of $\Seq{X_k}$ that has the same
distribution as $a^{-1} \rad \lit \lincr_1$.
Since $\lincr_k$ has the same law as $\lincr_1$ and is independent of
$X_{k-1}$, we have
$$
\sum_{k \ge 1} \Pbig{\rad \lit \lincr_k > a \cdot |X_{k-1}|}
=
\sum_{k \ge 1} \Pbig{R > |X_{k-1}|}
=
\EBig{\sum_{k \ge 1} \III{R > |X_{k-1}|}}
\,.
$$
The idea now is that for $k > R^2$, we control the chance that $|X_{k-1}| <
R$ by using \ref l.close-polyd/, summing over the relevant possible values of
$X_{k-1}$. Thus,
$$\eqaln{
\EBig{\sum_{k \ge 1} \III{R > |X_{k-1}|}}
&\le
\EBig{R^2 + \Ebig{\sum_{k > R^2} \III{R > |X_{k-1}|} \bigm| R}}
\cr&\le
\EBig{R^2 + \sum_{k > R^2} V_\baseg(R) c k^{-d/2}}
\cr&\le
\EBig{R^2 + c R^d (R^2)^{1-d/2}}
=
c \E[R^2]
\,.
\Qed\cr
}$$

\proofof t.2ndmoment
Our assumption is that $\Ebig{|\widehat X_1|^2} < \infty$.
However, all we will
use of this moment condition is weaker, namely,
that $H(\widehat X_1) < \infty$,
that $\Ebig{|X_1|^2} < \infty$,
and that $\Ebig{( \rad \lit \Lamps_1)^2} < \infty$.
The first is a well-known consequence of the weaker assumption
$\Ebig{|\widehat X_1|} < \infty$; the latter two
follow from $|(\phi, x)| \ge \max \big\{|x|, \rad \lit \phi\big\}
+ |\!\lit \phi|$.
\rl{Thus, we need not assume that $\lampg$ is finitely generated.}

The case $\E[X_1] \ne \bfz$ was done by \ref b.Kaim:survey/, so assume that
$\E[X_1] = \bfz$.

Let $\epsilon \in (0, 1)$.
Choose $s$ so large that
$$
\inf_n \P_\bfz \big[\all {k \le n} \|X_k\| \le s \sqrt n\,\big]
>
{1 \over 2} + {\Ebig{( \rad \lit \Lamps_1)^2} \over s^2}
\,;
$$
such an $s$ exists by \ref l.polylamps/.
We will define random sets $Q_{n, \epsilon}$ that are
$\Lamps_\infty$-measurable in order to apply \ref c.Kcrit-epsilon/.

Abbreviate $\baseg := \Z^d$.

Write $\ranen U := \Seq{\widehat Y_m \st n < m \le n + \epsilon n}$.
By \ref l.aep/,
there is a set $\enum U_n \subseteq \big(\lampg \wr \baseg \big)^{\flr{\epsilon
n}}$
with $\log |\enum U_n| < 2 \epsilon n H(\widehat X_1)$ and
$\P[\ranen U \in \enum U_n] \to 1$.

We wish to define a set
$Q_{n, \epsilon}(\Lamps_\infty)$ that will contain $\widehat X_n$ with
reasonable probability and that will have small exponential growth.
We will consider the possible increments $\ranen U$ and the possible
values of $X_n \in B(s \sqrt n\,)$.
Given such possible values, we guess values for $\Lamps_n$ from
the ones we see, $\Lamps_\infty$, by correcting by the changes caused by
$\ranen U$.
Namely, the lamp at some $z \in B(2s\sqrt n\,)$ is changed at time $m \in (n,
n+\epsilon n]$ by multiplying by $\lincr_m(z - X_{m-1})$, whence the
total change from what it was at time $n$ due to the changes in $\ranen U$ is
$\prod_{m=n+1}^{n + \flr{\epsilon n}} \lincr_m(z - X_{m-1})$. Provided
that the lamps in $B(2s\sqrt n\,)$ are not changed after time $n + \epsilon
n$, we may multiply $\Lamps_\infty(z)$ by the inverse of this product to
guess $\Lamps_n(z)$.

Thus, we proceed as follows.
Let $\lamps_\infty \in \lampg^{\baseg}$.
Define $Q_{n, \epsilon}(\lamps_\infty)$ to be the set of all
$(\lamps_n, x)$ such that
there is some
$U = \Seq{(\psi_{n+1}, y_{n+1}), \ldots, (\psi_{n+\flr{\epsilon n}}, y_{n +
\flr{\epsilon n}})}
\in \enum U_n$ and
some
$x \in B(s\sqrt n\,)$,
such that,
writing $z_m := x + \sum_{j=n+1}^m y_j$ for $n \le m \le n + \epsilon n$
and
$$
\psi(z) :=
\prod_{m=n+1}^{n + \flr{\epsilon n}} \psi_m(z - z_{m-1})
\,,
$$
we have
$$
\lamps_n(z) = \cases{
               \lamps_\infty(z)\psi(z)^{-1} & for $z \in B(2s\sqrt n\,)$,\cr
               \id & for $z \notin B(2s\sqrt n\,)$.\cr}
$$
By assumption, the number of choices of $U \in \enum U_n$ is at most $e^{
2 \epsilon n H(\widehat X_1)}$.
The number of choices of $x$ is at most $c n^{d/2}$.
Therefore, $|Q_{n, \epsilon}(\lamps_\infty)| < e^{2 \epsilon n H(\widehat
X_1) + o(n)}$.

We will prove that $\limsup_{n \to\infty} \Pbig{\widehat X_n \in
Q_{n, \epsilon}(\Lamps_\infty)} > 0$.


Let $A_n$ be the event that
$\|X_k\| \le s \sqrt n$ for all $k \le n$.
Let $C_n$ be the event that
$\rad \lit \lincr_k > s \sqrt n$ for some $k \le n$.
Then
$$
\P[\rad \lit \lincr_k > s \sqrt n]
\le
{\Ebig{( \rad \lit \lincr_k)^2} \over s^2 n}
=
{\Ebig{( \rad \lit \Lamps_1)^2} \over s^2 n}
$$
by Chebyshev's inequality, whence
$\P(C_n) \le \Ebig{( \rad \lit \Lamps_1)^2}/s^2$ by a union bound.
Let $D_n$ be the event that
$\Lamps_n(y) = \id$ for all $y \notin B(2s \sqrt n\,)$.
Then $A_n \setminus D_n \subseteq C_n$,
whence $\P(A_n \cap D_n) = \P(A_n) - \P(A_n \setminus D_n)
\ge \P(A_n) - \P(C_n) > 1/2$ by choice of $s$.

Let $E_n$ be the event that $\|X_{n + \flr{\epsilon n}} - X_n\| > 5s\sqrt n$.
Then $\liminf_{n \to\infty} \P(E_n) > 0$ and $E_n$ is independent of $A_n
\cap D_n$; on the event $A_n \cap D_n \cap E_n$, we have $\|X_{n +
\flr{\epsilon n}}\| > 4s\sqrt n$.
Let $F_n$ be the event that for all $m > n + \epsilon n$, we have
$\|X_m\| > 2s \sqrt n$.
By \ref l.polylamps/, $\lim_{n \to\infty} \P(F_n \mid A_nD_nE_n) =
1 - 1/2^{d-2}$.

Let $G_n$ be the event that at no time after $n + \flr{\epsilon n}$ does
the walk change a lamp in $B(2s\sqrt n\,)$.
Then $A_nD_nE_nF_n \setminus G_n$ is contained in the event that for some $m >
n+\epsilon n$, we have $\rad \lit \lincr_m > \|X_{m-1}\|/2$, which by the
Borel--Cantelli lemma and \ref l.separate/, has probability tending to 0 as
$n \to\infty$.
Therefore, $\liminf_{n \to\infty} \P(A_nD_nE_nF_n G_n) > 0$.

On the event $A_nD_nE_nF_n G_n$, we have that for every $z \in B(2s\sqrt
n\,)$,
$$
\Lamps_\infty(z)
=
\Lamps_n(z) \prod_{m=n+1}^{n + \flr{\epsilon n}} \lincr_m(z - X_{m-1})
\,,
$$
as desired.
\Qed

Recall that our proof of \ref t.2ndmoment/ did not use the full strength of
the hypothesis
$\Ebig{|\widehat X_1|^2} < \infty$, but only the weaker hypotheses
that $H(\widehat X_1) < \infty$,
that $\Ebig{|X_1|^2} < \infty$,
and that $\Ebig{( \rad \lit \Lamps_1)^2} < \infty$.
This last assumption cannot be weakened to finiteness of a smaller moment,
even if $\Seq{X_n}$ is simple random walk on $\Z^3$ and $\lampg = \Z_2$.
To see this, we adapt \ref b.Kaim:examples/, Proposition
1.1, which gave an example of a $\mu$-walk on $\Z_2 \wr \Z$ that yielded a
nontrivial Poisson boundary but with no limiting configuration of lamps
a.s.
Indeed, suppose that $\widehat X_1$ has the following distribution:
With probability 1/2, $\lincr_1 = \bfz$ and $X_1$ is a step of simple
random walk on $\Z^3$, while
for each $n \ge 1$, with probability $c_0/n^3$, $\, \lincr_1 = \I{B(n)}$ and
$X_1 = \bfz$,
where $c_0 := 1/(2\zeta(3))$ is a normalizing constant.
We still have $H(\widehat X_1) < \infty$, while
$\Ebig{( \rad \lit \Lamps_1)^a} < \infty$ iff $a < 2$.
We claim that while $\Lamps_\infty$ does not exist a.s.\ for this walk, the
Poisson boundary is nontrivial.
To see this,
condition on the walk in the base, $\Seq{X_n}$.
If $X_n = X_{n+1}$, then the chance that at time $n+1$
the lamp changes
at the origin is of order $1/\big(1 + \|X_n\|^2\big)$, independently
of all other steps of the walk.
Now $\sum_n \big(1 + \|X_n\|^2\big)^{-1} = \infty$ a.s.\ by the law of the
iterated logarithm,
whence the Borel--Cantelli lemma yields
infinitely many changes of the lamp at the origin a.s.
On the other hand, the difference between the lamp at the origin and the
lamp at $(1, 0, 0)$ changes only finitely many times a.s., again by the
Borel--Cantelli lemma, since
if $X_n = X_{n+1}$, then the chance that at time $n+1$ this difference
changes is of order $1/\big(1 + \|X_n\|^3\big)$, independently of all other
steps of the walk, and
$\sum_n \big(1 + \|X_n\|^3\big)^{-1} < \infty$ a.s.\ by \ref b.DvEr/.
Therefore, the Poisson boundary is nontrivial.
On the other hand, if $\mu$ has a finite first moment and projects to a
transient random walk on $\Z^d$, then a limiting lamp configuration exists;
see Theorem 3.3 of \ref b.Kaim:solvable/ or Lemma 1.1 of \ref b.Erschler/.
This general case is still open: is the harmonic measure on the limiting
lamp configuration equal to the Poisson boundary?
We remark that
\ref b.Erschler/ shows that the Poisson boundary can be nontrivial
even for some random walks where no combination of lamps stabilizes.

\bsection{Metabelian Groups}{s.meta}

As \ref b.Erschler/ noted following \ref b.Vershik:survey/, free metabelian
groups are sufficiently similar to lamplighter groups on $\Z^d$ that
similar results on their Poisson boundaries carry over. A group $F$ is
\dfn{metabelian} if $F''$ is trivial,
where prime indicates commutator subgroup.
Those of the form $\FF_d/\FF''_d$ are called \dfn{free metabelian groups},
where $\FF_d$ is the free group on $d$ generators.
More generally, consider groups of the form $\FF_d/H'$, where $H$ is a
normal subgroup of $\FF_d$.
As explained by \refbmulti{Erschler:Liouville}, with more details
given by \ref b.VershikDobrynin/, the groups $\FF_d/H'$ are isomorphic to
groups of finite configurations on $\baseg := \FF_d/H$ as follows.

Let $G$ be the right Cayley graph of $\FF_d/H$ corresponding to the free
generators of $\FF_d$.
Orient each edge of $G$ so as to form the group $C_1(G) = C_1(G, \Z)$ of
1-chains.
For each $x \in \baseg$, fix a finite path $\Seq{e_1, \ldots, e_k}$ of edges
from $o \in \baseg$ to $x$. To this path associate the 1-chain $\theta_x :=
\sum_{j=1}^k \pm e_j$, where we choose the plus sign iff $e_j$ is oriented
in the direction from $o$ to $x$ along the path.
For simplicity, we choose $\theta_o := 0$.
Let $Z_1(G)$ denote the space of cycles in $C_1(G)$.
(As there are no 2-cells, this is the same as $H_1(G, \Z)$.)
Note that $H$ is the fundamental group of $G$, and its abelianization,
$H/H'$, is canonically isomorphic to $Z_1(G)$, meaning that the
homomorphism $\varphi \colon \FF_d \to C_1(G)$ defined by
$\varphi(a) := \theta_{aH}$ for generators $a$ of $\FF_d$ has kernel $H'$ and
$\varphi(H) = Z_1(G)$.
Now $\baseg$ acts on $G$ by translation from the left, and so also acts on
$C_1(G)$, which we denote by $(x, f) \mapsto T_x f$.
Define $\widetilde \baseg$ to be the subset $\{\theta_x + f \st x \in
\baseg,\ f
\in Z_1(G)\} \subset C_1(G)$; this set is clearly
independent of the choices of the chains $\theta_x$.
In addition, the map $\theta_x + f \mapsto x$ from $\widetilde \baseg \to
\baseg$
is well defined.
Define a multiplication on $\widetilde \baseg$ by
$$
(\theta_x + f) (\theta_y + g)
:=
\theta_{x} + T_x \theta_y + f + T_x g
\,.
$$
Then $\widetilde \baseg$ is closed under this multiplication
because $\theta_x + T_x \theta_y$ corresponds to a path from $o$ to $xy$.
It is easy to check that $\widetilde \baseg$ is a group with
identity element $0$.
Indeed, $\widetilde \baseg$ is canonically isomorphic to $\FF_d/H'$ via
the homomorphism $\varphi$ defined above.

A random walk $\Seq{\theta_{X_n} + \Phi_n}$
on $\widetilde \baseg$
yields a.s.\ an edgewise limiting configuration
in the
space of cochains, $C^1(G)$, under weak conditions:
As \ref b.Erschler/ proved, it suffices that
the walk on $\widetilde \baseg$ has finite first moment and projects to a
transient random walk on $\baseg$.
Under similar conditions as our previous
theorems and with similar proofs, the subset of possible limits, together
with harmonic measures, is the Poisson boundary.
For example, if $\baseg = \FF_d/H$ has at least cubic growth, then this holds
for every finitely supported walk on $\widetilde \baseg$.
In the case of free metabelian groups with $d \ge 3$, it holds
for every walk having finite second moment.
\ref b.Erschler/ had proved this for free metabelian groups with $d \ge 5$
and $\mu$ having finite third moment. 

\medbreak
\noindent {\bf Acknowledgements.}\enspace
We are grateful to Vadim Kaimanovich and Anatoly
Vershik for posing the problem that
we address here and for bringing attention to it over the years.
We thank Anna Erschler for explaining her proof to us, Ori
Gurel-Gurevich for helpful discussions, and Omer Tamuz for a
useful reference.
We thank the referees for careful readings and expert
suggestions that led to
the improvement of our exposition.
\rl{We owe a particular debt to one of the referees 
for having urged us to extend our results.}

\def\noop#1{\relax}
\def\polhk#1{\setbox0=\hbox{#1}{\ooalign{\hidewidth
  \lower1.5ex\hbox{`}\hidewidth\crcr\unhbox0}}}
  \def\soft#1{\leavevmode\setbox0=\hbox{h}\dimen7=\ht0\advance \dimen7
  by-1ex\relax\if t#1\relax\rlap{\raise.6\dimen7
  \hbox{\kern.3ex\char'47}}#1\relax\else\if T#1\relax
  \rlap{\raise.5\dimen7\hbox{\kern1.3ex\char'47}}#1\relax \else\if
  d#1\relax\rlap{\raise.5\dimen7\hbox{\kern.9ex \char'47}}#1\relax\else\if
  D#1\relax\rlap{\raise.5\dimen7 \hbox{\kern1.4ex\char'47}}#1\relax\else\if
  l#1\relax \rlap{\raise.5\dimen7\hbox{\kern.4ex\char'47}}#1\relax \else\if
  L#1\relax\rlap{\raise.5\dimen7\hbox{\kern.7ex
  \char'47}}#1\relax\else\message{accent \string\soft \space #1 not
  defined!}#1\relax\fi\fi\fi\fi\fi\fi} 
  \def\lfhook#1{\setbox0=\hbox{#1}{\ooalign{\hidewidth
  \lower1.5ex\hbox{'}\hidewidth\crcr\unhbox0}}} 
   
\def\temp{\let\linkit=\linkyear \apaliketrue}
\temp
\ifcitationgeneration\immediate\write\labelfile{\sanitize\temp}\fi
\def\startreferences{
 \vskip0pt plus.3\vsize \penalty -150 \vskip0pt
 plus-.3\vsize \bigskip\bigskip \vskip \parskip
 \begingroup\baselineskip=12pt\frenchspacing
 \bibliographytitle
 \vskip12pt\parindent=0pt
 \def\and{{\rm and}}
 \def\em{\it}
 \def\newblock{\hskip .11em plus.33em minus.07em}
 \def\bibauthor##1{{\sc ##1}}
 \def\bibitem[##1]##2
 {\htmlanchor{##2}{}\RefLabel{##2}[##1]\hangindent=.8cm\hangafter=1}
 }
\def\endreferences{\bigskip\bigskip\endgroup}
\ifundefined{bibstylemodification}\relax\else\bibstylemodification\fi
\startreferences

\bibitem[Avez (1972)]{MR0324741}
\bibauthor{Avez, A.} (1972).
\newblock Entropie des groupes de type fini.
\newblock {\em C. R. Acad. Sci. Paris S\'er. A-B}, {\bf 275}, A1363--A1366.

\bibitem[Avez (1974)]{MR0353405}
\bibauthor{Avez, A.} (1974).
\newblock Th\'eor\`eme de {C}hoquet-{D}eny pour les groupes \`a croissance non
  exponentielle.
\newblock {\em C. R. Acad. Sci. Paris S\'er. A}, {\bf 279}, 25--28.

\bibitem[Avez (1976a)]{MR0482911}
\bibauthor{Avez, A.} (1976a).
\newblock Croissance des groupes de type fini et fonctions harmoniques.
\newblock In {\em Th\'eorie Ergodique}, Lecture Notes in Mathematics, Vol. 532,
  pages 35--49. Springer, Berlin.
\newblock Actes des Journ{\'e}es Ergodiques, Rennes, 1973/1974, Edit{\'e} par
  J.-P. Conze et M. S. Keane.

\bibitem[Avez (1976b)]{MR0507229}
\bibauthor{Avez, A.} (1976b).
\newblock Harmonic functions on groups.
\newblock In {\em Differential Geometry and Relativity}, pages 27--32.
  Mathematical Phys. and Appl. Math., Vol. 3. Reidel, Dordrecht.

\bibitem[Ballmann and Ledrappier (1994)]{MR1269841}
\bibauthor{Ballmann, W. \and{} Ledrappier, F.} (1994).
\newblock The {P}oisson boundary for rank one manifolds and their cocompact
  lattices.
\newblock {\em Forum Math.}, {\bf 6}(3), 301--313.

\bibitem[Benjamini and Peres (1994)]{MR94m:60141}
\bibauthor{Benjamini, I. \and{} Peres, Y.} (1994).
\newblock Tree-indexed random walks on groups and first passage percolation.
\newblock {\em Probab. Theory Related Fields}, {\bf 98}(1), 91--112.

\bibitem[Blach{\`e}re, Ha{\"{\i}}ssinsky, and Mathieu (2008)]{MR2408585}
\bibauthor{Blach{\`e}re, S., Ha{\"{\i}}ssinsky, P., \and{} Mathieu, P.} (2008).
\newblock Asymptotic entropy and {G}reen speed for random walks on countable
  groups.
\newblock {\em Ann. Probab.}, {\bf 36}(3), 1134--1152.

\bibitem[Blackwell (1955)]{MR17:754d}
\bibauthor{Blackwell, D.} (1955).
\newblock On transient {M}arkov processes with a countable number of states and
  stationary transition probabilities.
\newblock {\em Ann. Math. Statist.}, {\bf 26}, 654--658.

\bibitem[Brofferio and Schapira (2011)]{MR2837130}
\bibauthor{Brofferio, S. \and{} Schapira, B.} (2011).
\newblock Poisson boundary of {${\rm GL}_d(\Bbb Q)$}.
\newblock {\em Israel J. Math.}, {\bf 185}, 125--140.

\bibitem[Cartwright, Kaimanovich, and Woess (1994)]{MR1306556}
\bibauthor{Cartwright, D.I., Kaimanovich, V.A., \and{} Woess, W.} (1994).
\newblock Random walks on the affine group of local fields and of homogeneous
  trees.
\newblock {\em Ann. Inst. Fourier (Grenoble)}, {\bf 44}(4), 1243--1288.

\bibitem[Coulhon, Grigor'yan, and Pittet (2001)]{MR1871289}
\bibauthor{Coulhon, T., Grigor'yan, A., \and{} Pittet, C.} (2001).
\newblock A geometric approach to on-diagonal heat kernel lower bounds on
  groups.
\newblock {\em Ann. Inst. Fourier (Grenoble)}, {\bf 51}(6), 1763--1827.

\bibitem[Coulhon and Saloff-Coste (1993)]{MR94g:58263}
\bibauthor{Coulhon, T. \and{} Saloff-Coste, L.} (1993).
\newblock Isop\'erim\'etrie pour les groupes et les vari\'et\'es.
\newblock {\em Rev. Mat. Iberoamericana}, {\bf 9}(2), 293--314.

\bibitem[Cover and Thomas (2006)]{MR2239987}
\bibauthor{Cover, T.M. \and{} Thomas, J.A.} (2006).
\newblock {\em Elements of Information Theory}.
\newblock Wiley-Interscience, Hoboken, NJ, second edition.

\bibitem[Derriennic (1976)]{MR0423532}
\bibauthor{Derriennic, Y.} (1976).
\newblock Lois ``z\'ero ou deux'' pour les processus de {M}arkov.
  {A}pplications aux marches al\'eatoires.
\newblock {\em Ann. Inst. H. Poincar\'e Sect. B (N.S.)}, {\bf 12}(2), 111--129.

\bibitem[Derriennic (1980)]{MR588163}
\bibauthor{Derriennic, Y.} (1980).
\newblock Quelques applications du th\'eor\`eme ergodique sous-additif.
\newblock In {\em Journ\'ees sur les {M}arches {A}l\'eatoires}, volume 74 of
  {\em Ast\'erisque}, pages 183--201, 4. Soc. Math. France, Paris.
\newblock Held at Kleebach, March 5--10, 1979.

\bibitem[Doob (1959)]{MR0107098}
\bibauthor{Doob, J.L.} (1959).
\newblock Discrete potential theory and boundaries.
\newblock {\em J. Math. Mech.}, {\bf 8}, 433--458; erratum 993.

\bibitem[Durrett (2010)]{MR2722836}
\bibauthor{Durrett, R.} (2010).
\newblock {\em Probability: Theory and Examples}.
\newblock Cambridge Series in Statistical and Probabilistic Mathematics.
  Cambridge University Press, Cambridge, fourth edition.

\bibitem[Dvoretzky and Erd{\"o}s (1951)]{MR0047272}
\bibauthor{Dvoretzky, A. \and{} Erd{\"o}s, P.} (1951).
\newblock Some problems on random walk in space.
\newblock In {\em Proc. {S}econd {B}erkeley {S}ymposium on {M}ath. {S}tatist.
  and {P}robability, 1950}, pages 353--367. University of California Press,
  Berkeley.

\bibitem[Dynkin and Maljutov (1961)]{MR24:A1751}
\bibauthor{Dynkin, E.B. \and{} Maljutov, M.B.} (1961).
\newblock Random walk on groups with a finite number of generators.
\newblock {\em Dokl. Akad. Nauk SSSR}, {\bf 137}, 1042--1045.
\newblock English translation: {\it Soviet Math. Dokl.} (1961) {\bf 2},
  399--402.

\bibitem[Erschler (2004a)]{MR2144977}
\bibauthor{Erschler, A.} (2004a).
\newblock Boundary behavior for groups of subexponential growth.
\newblock {\em Ann. of Math. (2)}, {\bf 160}(3), 1183--1210.

\bibitem[Erschler (2004b)]{MR2025301}
\bibauthor{Erschler, A.} (2004b).
\newblock Liouville property for groups and manifolds.
\newblock {\em Invent. Math.}, {\bf 155}(1), 55--80.

\bibitem[Erschler (2010)]{MR2827814}
\bibauthor{Erschler, A.} (2010).
\newblock Poisson-{F}urstenberg boundaries, large-scale geometry and growth of
  groups.
\newblock In {\em Proceedings of the {I}nternational {C}ongress of
  {M}athematicians. {V}olume {II}}, pages 681--704. Hindustan Book Agency, New
  Delhi.

\bibitem[Erschler (2011)]{MR2745278}
\bibauthor{Erschler, A.} (2011).
\newblock Poisson-{F}urstenberg boundary of random walks on wreath products and
  free metabelian groups.
\newblock {\em Comment. Math. Helv.}, {\bf 86}(1), 113--143.

\bibitem[Feldman (1962)]{MR0139202}
\bibauthor{Feldman, J.} (1962).
\newblock Feller and {M}artin boundaries for countable sets.
\newblock {\em Illinois J. Math.}, {\bf 6}, 357--366.

\bibitem[Feller (1956)]{MR0090927}
\bibauthor{Feller, W.} (1956).
\newblock Boundaries induced by non-negative matrices.
\newblock {\em Trans. Amer. Math. Soc.}, {\bf 83}, 19--54.

\bibitem[Feller (1968)]{MR0228020}
\bibauthor{Feller, W.} (1968).
\newblock {\em An Introduction to Probability Theory and its Applications.
  {V}ol. {I}}.
\newblock Third edition. John Wiley \& Sons, Inc., New York-London-Sydney.

\bibitem[Forghani and Tiozzo (2019)]{MR4026614}
\bibauthor{Forghani, B. \and{} Tiozzo, G.} (2019).
\newblock Random walks of infinite moment on free semigroups.
\newblock {\em Probab. Theory Related Fields}, {\bf 175}(3--4), 1099--1122.

\bibitem[Frisch, Hartman, Tamuz, and Vahidi~Ferdowsi (2019)]{MR3990605}
\bibauthor{Frisch, J., Hartman, Y., Tamuz, O., \and{} Vahidi~Ferdowsi, P.}
  (2019).
\newblock Choquet-{D}eny groups and the infinite conjugacy class property.
\newblock {\em Ann. of Math. (2)}, {\bf 190}(1), 307--320.

\bibitem[Furstenberg (1963)]{MR0146298}
\bibauthor{Furstenberg, H.} (1963).
\newblock A {P}oisson formula for semi-simple {L}ie groups.
\newblock {\em Ann. of Math. (2)}, {\bf 77}, 335--386.

\bibitem[Furstenberg (1971a)]{MR0430160}
\bibauthor{Furstenberg, H.} (1971a).
\newblock Boundaries of {L}ie groups and discrete subgroups.
\newblock In {\em Actes du {C}ongr\`es {I}nternational des
  {M}ath\-\'e\-ma\-ti\-ciens ({N}ice, 1970), {T}ome 2}, pages 301--306.
  Gauthier-Villars, Paris.

\bibitem[Furstenberg (1971b)]{MR0284569}
\bibauthor{Furstenberg, H.} (1971b).
\newblock Random walks and discrete subgroups of {L}ie groups.
\newblock In {\em Advances in {P}robability and {R}elated {T}opics, {V}ol. 1},
  pages 1--63. Dekker, New York.

\bibitem[Furstenberg (1973)]{MR0352328}
\bibauthor{Furstenberg, H.} (1973).
\newblock Boundary theory and stochastic processes on homogeneous spaces.
\newblock In Moore, C.C., editor, {\em Harmonic Analysis on Homogeneous
  Spaces}, pages 193--229. Amer. Math. Soc., Providence, R.I.
\newblock Harmonic analysis on homogeneous spaces ({P}roc. {S}ympos. {P}ure
  {M}ath., {V}ol. {XXVI}, {W}illiams {C}oll., {W}illiamstown, {M}ass., 1972).

\bibitem[Gautero and Math{\'e}us (2012)]{MR3011489}
\bibauthor{Gautero, F. \and{} Math{\'e}us, F.} (2012).
\newblock Poisson boundary of groups acting on {$\Bbb R$}-trees.
\newblock {\em Israel J. Math.}, {\bf 191}(2), 585--646.

\bibitem[Georgakopoulos (2016)]{Agelos}
\bibauthor{Georgakopoulos, A.} (2016).
\newblock The boundary of a square tiling of a graph coincides with the
  {P}oisson boundary.
\newblock {\em Invent. Math.}, {\bf 203}(3), 773--821.

\bibitem[Gromov (1981)]{MR83b:53041}
\bibauthor{Gromov, M.} (1981).
\newblock Groups of polynomial growth and expanding maps.
\newblock {\em Inst. Hautes \'Etudes Sci. Publ. Math.}, {\bf 53}, 53--73.

\bibitem[Hunt (1960)]{MR0123364}
\bibauthor{Hunt, G.A.} (1960).
\newblock Markoff chains and {M}artin boundaries.
\newblock {\em Illinois J. Math.}, {\bf 4}, 313--340.

\bibitem[James and Peres (1996)]{MR1687097}
\bibauthor{James, N. \and{} Peres, Y.} (1996).
\newblock Cutpoints and exchangeable events for random walks.
\newblock {\em Teor. Veroyatnost. i Primenen.}, {\bf 41}(4), 854--868.

\bibitem[Kaimanovich (1983)]{MR697250}
\bibauthor{Kaimanovich, V.A.} (1983).
\newblock Examples of nonabelian discrete groups with nontrivial exit boundary.
\newblock {\em Zap. Nauchn. Sem. Leningrad. Otdel. Mat. Inst. Steklov. (LOMI)},
  {\bf 123}, 167--184.
\newblock Differential geometry, Lie groups and mechanics, V. English
  translation: {\it J. Soviet Math.} {\bf 28} (1985), no. 4, 579--591,
  \doi{10.1007/BF02104988}.

\bibitem[Kaimanovich (1985)]{MR780288}
\bibauthor{Kaimanovich, V.A.} (1985).
\newblock An entropy criterion of maximality for the boundary of random walks
  on discrete groups.
\newblock {\em Dokl. Akad. Nauk SSSR}, {\bf 280}(5), 1051--1054.

\bibitem[Kaimanovich (1991)]{MR1178986}
\bibauthor{Kaimanovich, V.A.} (1991).
\newblock Poisson boundaries of random walks on discrete solvable groups.
\newblock In Heyer, H., editor, {\em Probability Measures on Groups. {X}},
  pages 205--238. Plenum, New York.

\bibitem[Kaimanovich (1994)]{MR1260536}
\bibauthor{Kaimanovich, V.A.} (1994).
\newblock The {P}oisson boundary of hyperbolic groups.
\newblock {\em C. R. Acad. Sci. Paris S\'er. I Math.}, {\bf 318}(1), 59--64.

\bibitem[Kaimanovich (2000)]{MR1815698}
\bibauthor{Kaimanovich, V.A.} (2000).
\newblock The {P}oisson formula for groups with hyperbolic properties.
\newblock {\em Ann. of Math. (2)}, {\bf 152}(3), 659--692.

\bibitem[Kaimanovich (2001)]{Kaim:survey}
\bibauthor{Kaimanovich, V.A.} (2001).
\newblock Poisson boundary of discrete groups.
\newblock Preprint, \hfil\break
  \url{http://citeseerx.ist.psu.edu/viewdoc/summary?doi=10.1.1.6.6675}.

\bibitem[Kaimanovich and Masur (1996)]{MR1395719}
\bibauthor{Kaimanovich, V.A. \and{} Masur, H.} (1996).
\newblock The {P}oisson boundary of the mapping class group.
\newblock {\em Invent. Math.}, {\bf 125}(2), 221--264.

\bibitem[Kaimanovich and Masur (1998)]{MR1636940}
\bibauthor{Kaimanovich, V.A. \and{} Masur, H.} (1998).
\newblock The {P}oisson boundary of {T}eichm\"uller space.
\newblock {\em J. Funct. Anal.}, {\bf 156}(2), 301--332.

\bibitem[Kaimanovich and Vershik (1983)]{MR85d:60024}
\bibauthor{Kaimanovich, V.A. \and{} Vershik, A.M.} (1983).
\newblock Random walks on discrete groups: {B}oundary and entropy.
\newblock {\em Ann. Probab.}, {\bf 11}(3), 457--490.

\bibitem[Kaimanovich and Woess (2002)]{MR1894110}
\bibauthor{Kaimanovich, V.A. \and{} Woess, W.} (2002).
\newblock Boundary and entropy of space homogeneous {M}arkov chains.
\newblock {\em Ann. Probab.}, {\bf 30}(1), 323--363.

\bibitem[Karlsson (2003)]{MR2011923}
\bibauthor{Karlsson, A.} (2003).
\newblock Boundaries and random walks on finitely generated infinite groups.
\newblock {\em Ark. Mat.}, {\bf 41}(2), 295--306.

\bibitem[Karlsson and Ledrappier (2007)]{MR2402595}
\bibauthor{Karlsson, A. \and{} Ledrappier, F.} (2007).
\newblock Linear drift and {P}oisson boundary for random walks.
\newblock {\em Pure Appl. Math. Q.}, {\bf 3}(4, Special Issue: In honor of
  Grigory Margulis. Part 1), 1027--1036.

\bibitem[Karlsson and Woess (2007)]{MR2318539}
\bibauthor{Karlsson, A. \and{} Woess, W.} (2007).
\newblock The {P}oisson boundary of lamplighter random walks on trees.
\newblock {\em Geom. Dedicata}, {\bf 124}, 95--107.

\bibitem[Lawler and Limic (2010)]{MR2677157}
\bibauthor{Lawler, G.F. \and{} Limic, V.} (2010).
\newblock {\em Random Walk: A Modern Introduction}, volume 123 of {\em
  Cambridge Studies in Advanced Mathematics}.
\newblock Cambridge University Press, Cambridge.

\bibitem[Ledrappier (1983)]{MR699165}
\bibauthor{Ledrappier, F.} (1983).
\newblock Une relation entre entropie, dimension et exposant pour certaines
  marches al\'eatoires.
\newblock {\em C. R. Acad. Sci. Paris S\'er. I Math.}, {\bf 296}(8), 369--372.

\bibitem[Ledrappier (1985)]{MR800190}
\bibauthor{Ledrappier, F.} (1985).
\newblock Poisson boundaries of discrete groups of matrices.
\newblock {\em Israel J. Math.}, {\bf 50}(4), 319--336.

\bibitem[Lyons and Oveis~Gharan (2018)]{MR3892273}
\bibauthor{Lyons, R. \and{} Oveis~Gharan, S.} (2018).
\newblock Sharp bounds on random walk eigenvalues via spectral embedding.
\newblock {\em Int. Math. Res. Not. IMRN}, {\bf 2018}(24), 7555--7605.

\bibitem[Lyons and Peres (2016)]{LP:book}
\bibauthor{Lyons, R. \and{} Peres, Y.} (2016).
\newblock {\em Probability on Trees and Networks}, volume 42 of {\em Cambridge
  Series in Statistical and Probabilistic Mathematics}.
\newblock Cambridge University Press, New York.
\newblock Available at \url{http://pages.iu.edu/~rdlyons/}.

\bibitem[Maher and Tiozzo (2018)]{MR3849626}
\bibauthor{Maher, J. \and{} Tiozzo, G.} (2018).
\newblock Random walks on weakly hyperbolic groups.
\newblock {\em J. Reine Angew. Math.}, {\bf 742}, 187--239.

\bibitem[Malyutin, Nagnibeda, and Serbin (2017)]{MR3644015}
\bibauthor{Malyutin, A., Nagnibeda, T., \and{} Serbin, D.} (2017).
\newblock Boundaries of {$\Bbb Z^n$}-free groups.
\newblock In Ceccherini-Silberstein, T., Salvatori, M., \and{} Sava-Huss, E.,
  editors, {\em Groups, Graphs and Random Walks}, volume 436 of {\em London
  Math. Soc. Lecture Note Ser.}, pages 355--390. Cambridge Univ. Press,
  Cambridge.
\newblock Selected papers from the workshop held in Cortona, June 2--6, 2014.

\bibitem[Malyutin (2003)]{MR2032055}
\bibauthor{Malyutin, A.V.} (2003).
\newblock The {P}oisson-{F}urstenberg boundary of a locally free group.
\newblock {\em Zap. Nauchn. Sem. S.-Peterburg. Otdel. Mat. Inst. Steklov.
  (POMI)}, {\bf 301}(Teor. Predst. Din. Sist. Komb. i Algoritm. Metody. 9),
  195--211, 245.

\bibitem[Malyutin and Svetlov (2014)]{MalSvet}
\bibauthor{Malyutin, A.V. \and{} Svetlov, P.} (2014).
\newblock {P}oisson-{F}urstenberg boundaries of fundamental groups of closed
  3-manifolds.
\newblock Preprint, \arXiv{1403.2135}.

\bibitem[Morris and Peres (2005)]{MR2198701}
\bibauthor{Morris, B. \and{} Peres, Y.} (2005).
\newblock Evolving sets, mixing and heat kernel bounds.
\newblock {\em Probab. Theory Related Fields}, {\bf 133}(2), 245--266.

\bibitem[M{\"o}rters and Peres (2010)]{BMbook}
\bibauthor{M{\"o}rters, P. \and{} Peres, Y.} (2010).
\newblock {\em Brownian Motion}, volume 30 of {\em Cambridge Series in
  Statistical and Probabilistic Mathematics}.
\newblock Cambridge University Press, Cambridge.
\newblock With an appendix by Oded Schramm and Wendelin Werner.

\bibitem[Nevo and Sageev (2013)]{MR3095714}
\bibauthor{Nevo, A. \and{} Sageev, M.} (2013).
\newblock The {P}oisson boundary of {${\rm CAT}(0)$} cube complex groups.
\newblock {\em Groups Geom. Dyn.}, {\bf 7}(3), 653--695.

\bibitem[Rosenblatt (1981)]{MR630645}
\bibauthor{Rosenblatt, J.} (1981).
\newblock Ergodic and mixing random walks on locally compact groups.
\newblock {\em Math. Ann.}, {\bf 257}(1), 31--42.

\bibitem[Sava (2010a)]{Sava:thesis}
\bibauthor{Sava, E.} (2010a).
\newblock {\em Lamplighter Random Walks and Entropy-Sensitivity of Languages}.
\newblock Ph.D. thesis, Technische Universit\"at Graz.
\newblock Available at \arXiv{1012.2757}.

\bibitem[Sava (2010b)]{MR2600904}
\bibauthor{Sava, E.} (2010b).
\newblock A note on the {P}oisson boundary of lamplighter random walks.
\newblock {\em Monatsh. Math.}, {\bf 159}(4), 379--396.

\bibitem[Spitzer (1976)]{MR52:9383}
\bibauthor{Spitzer, F.} (1976).
\newblock {\em Principles of {R}andom {W}alk}, volume 34 of {\em Graduate Texts
  in Mathematics}.
\newblock Springer-Verlag, New York, second edition.

\bibitem[Varopoulos (1985)]{MR87j:60100}
\bibauthor{Varopoulos, N.{\relax Th}.} (1985).
\newblock Long range estimates for {M}arkov chains.
\newblock {\em Bull. Sci. Math. (2)}, {\bf 109}(3), 225--252.

\bibitem[Vershik (2000)]{MR1786730}
\bibauthor{Vershik, A.M.} (2000).
\newblock Dynamic theory of growth in groups: {E}ntropy, boundaries, examples.
\newblock {\em Uspekhi Mat. Nauk}, {\bf 55}(4(334)), 59--128.

\bibitem[Vershik and Dobrynin (2005)]{MR2197831}
\bibauthor{Vershik, A.M. \and{} Dobrynin, S.V.} (2005).
\newblock Geometrical approach to the free solvable groups.
\newblock {\em Internat. J. Algebra Comput.}, {\bf 15}(5--6), 1243--1260.

\bibitem[Vershik and Kaimanovich (1979)]{MR553972}
\bibauthor{Vershik, A.M. \and{} Kaimanovich, V.A.} (1979).
\newblock Random walks on groups: {B}oundary, entropy, uniform distribution.
\newblock {\em Dokl. Akad. Nauk SSSR}, {\bf 249}(1), 15--18.

\bibitem[Whitt (2002)]{MR1876437}
\bibauthor{Whitt, W.} (2002).
\newblock {\em Stochastic-Process Limits}.
\newblock Springer Series in Operations Research. Springer-Verlag, New York.
\newblock An introduction to stochastic-process limits and their application to
  queues.

\endreferences

\filbreak
\begingroup
\eightpoint\sc
\parindent=0pt\baselineskip=10pt

Department of Mathematics,
831 E.\ 3rd St.,
Indiana University,
Bloomington, IN 47405-7106
\emailwww{rdlyons@indiana.edu}
{http://pages.iu.edu/\string~rdlyons/}

\smallskip

Kent State University,
Department of Mathematical Sciences,
Mathematics and Computer Science Building 233,
Summit Street, 
Kent, OH 44242
\emailwww{yperes@gmail.com}
{https://yuvalperes.com/}

\endgroup

\bye